\documentclass[11pt]{article}

\usepackage[T1]{fontenc}
\usepackage[latin1]{inputenc}

\usepackage[]{graphicx,float,latexsym,times}
\usepackage{amsfonts,amstext,amsmath,amssymb,amsthm,dsfont}
 
\usepackage{enumerate}
\usepackage{setspace}

\RequirePackage{xcolor}
\RequirePackage[numbers]{natbib}
\RequirePackage[colorlinks,citecolor=blue,urlcolor=blue]{hyperref}

\usepackage{setspace}

\usepackage{epsfig}
\usepackage{graphicx}
\usepackage{wrapfig}
\usepackage{amssymb}
\usepackage{amsfonts}
\usepackage{amsmath}
\usepackage{amsthm}
\usepackage{enumerate}
\usepackage{multicol}
\usepackage{bbm} 

\newtheorem{theorem}{Theorem}[section]
\newtheorem{proposition}[theorem]{Proposition}
\newtheorem{lemma}[theorem]{Lemma}

\newtheorem{remark}{Remark}[section]
\newtheorem{corollary}[theorem]{Corollary}
\newcommand\fdem{\hfill $\Box$}
\newcommand\cA{{\cal A}}
\newcommand\cC{{\cal C}}

\newcommand\cG{{\cal G}}

\newcommand\cD{{\cal D}}
\newcommand\cQ{{\cal Q}}

\newcommand\cR{{\cal R}}

\newcommand\ve{\varepsilon}

\def\bbr{{\mathbb R}}

\def\bbc{{\mathbb C}}

\def\text#1{\hbox{#1}}
\def\proof{{\noindent \bf Proof. }}

\def\E{{\bf E}}

\def\P{{\bf P}}
\def\C{{\bf C}}
\def\c{{\bf c}}
\def\D{{\bf D}}
\def\H{{\bf H}}

\def\U{{\bf U}}

\def\r{{\bf r}}
\def\l{{\bf l}}

\def\d{{\bf d}}
\def\L{{\bf L}}

\newcommand{\wh}{\widehat}
\newcommand{\wt}{\widetilde}
\newcommand\oo{\mbox{o}}

\newcommand\Er{\mbox{Err}}

\newcommand\Tr{\mbox{Tr}}

\def\Chi{{\bf 1}}
\def\d{\mathrm{d}}
\def\build #1_#2{\mathrel{\mathop{\kern 0pt #1}\limits_\zs{#2}}}
\newcommand{\zs}[1]{{\mathchoice{#1}{#1}{\lower.25ex\hbox{$\scriptstyle#1$}}
{\lower0.25ex\hbox{$\scriptscriptstyle#1$}}}}
\numberwithin{equation}{section}

\begin{document}
\title{
Robust adaptive efficient  estimation for semi-Markov nonparametric regression
models
\thanks{
This work was done under partial financial support
of the grant of RSF number 14-49-00079 (National Research University ``MPEI''
14 Krasnokazarmennaya, 111250 Moscow, Russia)
and of the RFBR Grant 16-01-00121.
}
}

\author{Vlad Stefan Barbu\thanks{
Laboratoire de Math\'ematiques Rapha\"el Salem,
 UMR 6085 CNRS-Universit\'e de Rouen,  France, e-mail:
barbu@univ-rouen.fr
}
,
Slim Beltaief\thanks{
Laboratoire de Math\'ematiques Rapha\"el Salem,
 UMR 6085 CNRS-Universit\'e de Rouen,  France, e-mail:
slim.beltaief1@univ-rouen.fr}
 and 
 Serguei Pergamenshchikov\thanks{
 Laboratoire de Math\'ematiques Rapha\"el Salem,
 UMR 6085 CNRS-Universit\'e de Rouen,  France
and
International Laboratory of Statistics of Stochastic Processes and
Quantitative Finance of National Research Tomsk State University,
 e-mail:
Serge.Pergamenshchikov@univ-rouen.fr}
}
\date{}

\maketitle

\begin{abstract}
We consider the  nonparametric robust estimation problem for regression models
 in continuous time with  semi-Markov noises. An adaptive model selection
procedure  is proposed. Under general moment conditions on the noise distribution 
a sharp non-asymptotic oracle inequality for the robust risks is obtained
 and the robust efficiency  is shown. It turns out that for semi-Markov models  the
 robust   minimax convergence  rate may be faster or slower  than the classical one.
 \end{abstract}

\vspace*{5mm}
\noindent {\sl MSC:} primary 62G08, secondary 62G05

\vspace*{5mm}
\noindent {\sl Keywords}: Non-asymptotic estimation; Robust risk;
 Model selection; Sharp oracle
inequality; Asymptotic efficiency.

\newpage

\section{Introduction}\label{sec:In}

Let us consider a regression model in continuous time
\begin{equation}\label{sec:In.1}
 \d\,y_\zs{t} = S(t)\d\,t + \d\,\xi_\zs{t}\,,\quad 
 0\le t \le n\,,
\end{equation}
where $S(\cdot)$ is an unknown $1$-periodic function from $\L_\zs{2}[0,1]$ defined on $\bbr$ with values in $\bbr$,  
the noise process $(\xi_\zs{t})_\zs{t\ge\, 0}$ is defined  as 
\begin{equation}\label{sec:Ex.1}
 \xi_\zs{t} = \varrho_\zs{1} L_\zs{t} + \varrho_\zs{2} z_\zs{t}\,,
\end{equation}
where $\varrho_\zs{1}$ and $\varrho_\zs{2}$ are unknown coefficients,
$(L_\zs{t})_\zs{t\ge\,0}$ is a Levy process and the pure jump process $(z_\zs{t})_\zs{t\ge\,1},$ defined in \eqref{sec:Ex.2}, is assumed to be a semi-Markov process (see, for example, \cite{BarbuLimnios2008}).

The problem is to estimate the unknown function $S$ in the model
\eqref{sec:In.1} on the basis of observations $(y_\zs{t})_\zs{0\le t\le n}$. Firstly, this problem was considered
in the framework of the ``signal+white noise'' models (see, for example, \cite{IbragimovKhasminskii1981} or \cite{Pinsker1981}).
Later, in order to study dependent observations in continuous time, were introduced ``signal+color noise''  regressions based
 on Ornstein-Uhlenbeck processes (cf. \cite{HopfnerKutoyants2009}, \cite{HopfnerKutoyants2010}, \cite{KonevPergamenshchikov2003}, \cite{KonevPergamenshchikov2010}).

Moreover, to include jumps in such models,  the papers 
\cite{KonevPergamenshchikov2012} and \cite{KonevPergamenshchikov2015}
used non Gaussian Ornstein-Uhlenbeck processes introduced in
\cite{BarndorffNielsenShephard2001} for modeling
of the risky assets in the stochastic volatility
financial markets. Unfortunately, the dependence of the
stable Ornstein-Uhlenbeck type decreases with a geometric rate.
So, asymptotically when the duration of
observations
   goes to infinity, 
 we obtain
very quickly
 the same
 ``signal+white noise'' model.

The main goal of this paper is to consider continuous time
 regression models with  dependent observations for which the dependence does not disappear for a sufficient large
 duration of observations. To this end we define the noise in the model \eqref{sec:In.1}
through a semi-Markov process which keeps  the dependence for any duration $n$. This type of models allows, for example, 
to estimate the signals observed under long impulse noise impact with a memory or ``against signals''.

In this paper we use the robust estimation approach introduced in \cite{KonevPergamenshchikov2012}
for such problems. To this end, we denote by $Q$ the distribution of $(\xi_\zs{t})_\zs{0\le t\le n}$ in the Skorokhod space $\cD[0,n]$.
 We assume that  $Q$ is unknown and  belongs to  some distribution family $\cQ_\zs{n}$ specified  in Section \ref{sec:Mrs}.
 In this paper we use the  quadratic risk
\begin{equation}\label{sec:In.4}
\cR_\zs{Q}(\wt{S}_\zs{n},S)=
\E_\zs{Q,S}\,\|\wt{S}_\zs{n}-S\|^{2}\,,
\end{equation}
where 
 $\|f\|^{2}=\int^{1}_\zs{0}\,f^{2}(s)\d s$ and
$\E_\zs{Q,S}$ is the expectation with respect to the distribution $\P_\zs{Q,S}$
 of the process  \eqref{sec:In.1} corresponding to the noise distribution $Q$. 
Since the noise distribution $Q$ is unknown, it seems reasonable
to introduce the
robust risk of the form
\begin{equation}\label{sec:In.6}
\cR^{*}_\zs{n}(\wt{S}_\zs{n},S)=\sup_\zs{Q\in\cQ_\zs{n}}\,
\cR_\zs{Q}(\wt{S}_\zs{n},S)\,,
\end{equation}
which enables us to take into account the information that
 $Q\in\cQ_\zs{n}$ and ensures the quality of an estimate $\wt{S}_\zs{n}$ for all
 distributions in the family $\cQ_\zs{n}$.

To summarize, the goal of this paper is to develop robust efficient model selection methods
for the model \eqref{sec:In.1} with the semi-Markov noise having unknown distribution, based on the approach proposed by Konev and Pergamenshchikov in 
 \cite{KonevPergamenshchikov2012} 
  and
  \cite{KonevPergamenshchikov2015} 
 for continuous time regression models with semimartingale noises.
   Unfortunately, we cannot use directly this method  for semi-Markov regression models, since their tool essentially uses
  the fact that the Ornstein-Uhlenbeck dependence decreases with geometrical rate and the ``white noise'' case is obtained sufficiently quickly.  
	
	Thus in the present  paper we propose new analytical tools based on renewal methods
  to obtain the sharp non-asymptotic oracle inequalities. As a consequence, we obtain the robust efficiency for the proposed model selection procedures
  in the adaptive setting.
  
  The rest of the paper is organized as follows. We start by introducing the main conditions in the next section. Then, in Section \ref{sec:Mo} we
construct the model selection procedure on the basis of
the
weighted least squares estimates. The main results are stated in Section \ref{sec:Mrs}; here we also specify the set
 of admissible weight sequences in the model selection procedure. In Section \ref{sec:Rtl}
 we derive some renewal results useful for obtaining other results of the paper.  In Section  \ref{sec:Smp}
 we develop stochastic calculus  for semi Markov processes. In Section~\ref{sec:Prsm} we study some properties of the model \eqref{sec:In.1}. 
  A numerical example is presented in Section \ref{sec:Siml}. Most of the results of the paper are proved in Section \ref{sec:Pr}. In Appendix some auxiliary propositions are given.

\bigskip

\section{Main conditions}\label{sec:Mcs}

In the model \eqref{sec:Ex.1} we assume that the Levy  process $L_\zs{t}$ is defined as 
 \begin{equation}\label{sec:Mcs.1}
L_\zs{t}=\check{\varrho}\,w_\zs{t}
+
\sqrt{1-\check{\varrho}^{2}}\,\check{L}_\zs{t}\,,\quad
\check{L}_\zs{t}=x*(\mu-\wt{\mu})_\zs{t}
\,,
\end{equation}
 where, $0\le \check{\varrho}\le 1$ is an unknown constant,
$(w_\zs{t})_\zs{t\ge\,0}$ is a standard Brownian motion, 
 $\mu(\d s,\d x)$ is the jump measure with the deterministic
compensator $\wt{\mu}(\d s\,\d x)=\d s\Pi(\d x)$, where 
$\Pi(\cdot)$ is some positive measure on $\bbr$ 
(see, for example
\cite{JacodShiryaev2002, ContTankov2004} for details) for which we assume that
\begin{equation}\label{sec:Mcs.2}
\Pi(x^{2})
=1
\quad\mbox{and}\quad
\Pi(x^{8})
\,<\,\infty\,,
\end{equation}
where we use the usual notation  $\Pi(\vert x\vert^{m})=\int_\zs{\bbr}\,\vert z \vert^{m}\,\Pi(\d z)$ for any $m>0$.
Note that $\Pi(\bbr)$ may be equal to $+\infty$.
Moreover, we assume that the pure jump process $(z_\zs{t})_\zs{t\ge\, 0}$ in
 \eqref{sec:Ex.1}
 is a semi-Markov process with the following form
 \begin{equation}\label{sec:Ex.2}
 z_\zs{t} = \sum_\zs{i=1}^{N_\zs{t}} Y_\zs{i},
\end{equation}
where $(Y_\zs{i})_\zs{i\ge\, 1}$ is an i.i.d. sequence of random variables with
\begin{equation*}\label{sec:Ex.3}
\E\,Y_\zs{i}=0\,,\quad 
\E\,Y^2_\zs{i}=1
\quad\mbox{and}\quad
\E\,Y^4_\zs{i}<\infty
\,.
\end{equation*}
Here $N_\zs{t}$ is a general  counting process (see, for example, \cite{Mikosch2004})
defined as
\begin{equation}\label{sec:Ex.4}
N_\zs{t} = \sum_\zs{k=1}^{\infty} \mathbbm{1}_\zs{\{T_\zs{k} \le t\}}
\quad\mbox{and}\quad
T_\zs{k}=\sum_\zs{l=1}^k\, \tau_\zs{l}\,,
\end{equation}
where $(\tau_\zs{l})_\zs{l\ge\,1}$ is an i.i.d. sequence of positive integrated
 random variables with distribution $\eta$ and mean $\check{\tau}=\E\,\tau_\zs{1}>0$. We assume that the processes 
$(N_\zs{t})_\zs{t\ge 0}$ and  $(Y_\zs{i})_{i\ge\, 1}$ are independent between them and are also independent of $(L_\zs{t})_\zs{t\ge 0}$.

Note that the process $(z_\zs{t})_\zs{t\ge\, 0}$ is a special case of a semi-Markov process (see, e.g., \cite{BarbuLimnios2008} and \cite{LO}).

\begin{remark} \label{Re.sec.Ex.0}
 It should be noted that if $\tau_\zs{j}$ are exponential random variables, then $(N_\zs{t})_\zs{t\ge 0}$
is a Poisson process and, in this case, $(\xi_\zs{t})_\zs{t\ge 0}$ is a Levy process for which this model has been studied
in   \cite{KonevPergamenshchikov2009a}, \cite{KonevPergamenshchikov2009b} and \cite{KonevPergamenshchikov2012}.
But, in the general case when the process \eqref{sec:Ex.2} is not a Levy process, this process has a memory and cannot be treated
in the framework of semi-martingales with independent increments.  In this case, we need to 
develop new tools based on renewal theory arguments, what we do in Section \ref{sec:Rtl}. This tools will be intensively used in the proofs of the main results of this paper.  
\end{remark}

Note that for any function $f$ from $\L_\zs{2}[0,n],$ $f: [0,n] \to\bbr,$ 
for the noise process $(\xi_\zs{t})_\zs{t\ge\, 0}$ defined in \eqref{sec:Ex.1}, with $(z_\zs{t})_\zs{t\ge\, 0}$ given in
 \eqref{sec:Ex.2}, the integral
\begin{equation}\label{sec:In.2}
I_\zs{n}(f)=\int_\zs{0}^{n} f(s) \d\xi_\zs{s}
\end{equation}
is well defined with $\E_\zs{Q}\,I_n(f)=0$. Moreover, as it is shown in Lemma \ref{Le.sec:Smp.1},
\begin{equation}\label{sec:In.3}
\E_\zs{Q}\,I^{2}_n(f) \le \varkappa_\zs{Q}\,\int_\zs{0}^{n} f^2(s) \d\,s\,,
\end{equation}
where $\varkappa_\zs{Q}=\varrho^{2}_\zs{1}+\varrho^{2}_\zs{2}\,\vert\rho\vert_\zs{*}$ and 
$\vert\rho\vert_\zs{*}=\sup_\zs{t\ge 0}\vert\rho(t)\vert<\infty$. Here $\rho$ is the density of
 the renewal measure 
$\check{\eta}$ defined as
\begin{equation}\label{sec:Cns.1}
\check{\eta}
=\sum^{\infty}_\zs{l=1}\,\eta^{(l)}
\,,
\end{equation}
where $\eta^{(l)}$ is the $l$th convolution power for  $\eta$. 
To study the series  \eqref{sec:Cns.1} we assume that the measure  $\eta$  has a density $g$ which satisfies the following conditions.

\bigskip

$\H_\zs{1}$) {\em Assume that, for any $x\in\bbr,$ 
there exist the finite limits
$$
g(x-)=\lim_\zs{z\to x-}g(z)
\quad\mbox{and}\quad
g(x+)=\lim_\zs{z\to x+}g(z)
$$
and, for any $K>0,$
there exists $\delta=\delta(K)>0$
for which
\begin{equation*}
\sup_\zs{\vert x\vert\le K}\,
\int^{\delta}_\zs{0}\,
\frac{
\vert 
g(x+t)+g(x-t)-g(x+)-g(x-)
\vert
}{t}
\d t
\,<\,\infty.
\end{equation*}
}

\bigskip

$\H_\zs{2}$) {\em  For any $\gamma>0,$
$$
\sup_\zs{z\ge 0}\,z^{\gamma}\vert 2g(z) -g(z-)-g(z+) \vert\,<\,\infty.
$$
}

$\H_\zs{3}$) {\em  There exists $\beta>0$ such that $\int_\zs{\bbr}\,e^{\beta x}\,g(x)\,\d x<\infty.$
}

\bigskip

\begin{remark} \label{Re.sec.Ex.2}
It should be noted that the condition $\H_\zs{3})$ means that there exists an exponential moment for the random variable $(\tau_\zs{j})_\zs{j\ge 1}$, i.e. these random variables are not too large.
 This is a natural 
constraint since these random variables define the intervals between jumps, i.e., the frequency of the jumps. 
 So, to study the influence of the jumps in the model \eqref{sec:In.1} one needs to consider the noise process \eqref{sec:Ex.1} with ``small'' interval between jumps or
 large jump frequency. 
\end{remark}

For the next condition we need to introduce the Fourier transform of any function $f$ from $\L_\zs{1}(\bbr),$ $f : \bbr\to\bbr,$  defined 
 as
\begin{equation}\label{sec:Rtl.06-0}
\wh{f}(\theta)=\frac{1}{2\pi}\,\int_\zs{\bbr}\,e^{i\theta x}\,f(x)\,\d x.
\end{equation}

\bigskip

$\H_\zs{4}$) {\em  There exists $t^{*}>0$ such that  the function $\wh{g}(\theta-it)$ 
belongs to $\L_\zs{1}(\bbr)$ for any $0\le t\le t^{*}$.
}

It is clear that Conditions $\H_\zs{1})$--$\H_\zs{4})$ hold true for any continuously differentiable function $g$, for example for the exponential density.
 
Now we define the family of the noise distributions
for the model \eqref{sec:In.1} which is
  used in the robust risk \eqref{sec:In.6}. Note that any distribution $Q$ from $\cQ_\zs{n}$
is defined by the unknown parameters  in  \eqref{sec:Ex.1} and \eqref{sec:Mcs.1}. 
We assume that 
\begin{equation}\label{sec:Ex.5}
\varsigma_\zs{*}\le \varrho^{2}_\zs{1}\le\varsigma^{*}\,,\quad 0\le \check{\varrho}\le 1
\quad\mbox{and}\quad
\varsigma_\zs{*}\le
\sigma_\zs{Q}
\le \varsigma^{*},
\end{equation}
where $\sigma_\zs{Q}=\varrho^{2}_\zs{1}+ \varrho^{2}_\zs{2}/\check{\tau}$, the    unknown 
bounds $0<\varsigma_\zs{*}\le \varsigma^{*}$ are functions of $n$, i.e.  $\varsigma_\zs{*}=\varsigma_\zs{*}(n)$ and  $\varsigma^{*}=\varsigma^{*}(n)$, such that
for any $\check{\epsilon}>0,$
\begin{equation}\label{sec:Mrs.5-1}
\lim_\zs{n\to\infty}n^{\check{\epsilon}}\,\varsigma_\zs{*}(n)=+\infty
\quad\mbox{and}\quad
\lim_\zs{n\to\infty}\,\frac{\varsigma^{*}(n)}{n^{\check{\epsilon}}}=0\,.
\end{equation}

\begin{remark} \label{Re.sec.Ex.1}
As we will see later,  the parameter $\sigma_\zs{Q}$ is the limit for the Fourier transform of the noise process \eqref{sec:Ex.1}.
 Such limit is called variance proxy
 (see \cite{KonevPergamenshchikov2012}).
\end{remark}

\begin{remark} \label{Re.sec.Ex.22++}
Note that, generally (but it is not necessary) the parameters $\varrho_\zs{1}$ and $\varrho_\zs{2}$ can be dependent on $n$. The conditions
\eqref{sec:Mrs.5-1} means that we consider all possible cases, i.e. these parameters  may go to the infinity or be constant or to the zero as well. 
See, for example, the conditions
(3.32) in
\cite{KonevPergamenshchikov2015}.
\end{remark}

\section{Model selection}\label{sec:Mo}

\noindent 
Let $(\phi_\zs{j})_\zs{j\ge\, 1}$ be an orthonormal uniformly bounded basis in $\L_\zs{2}[0,1]$, i.e., 
for some  constant $\phi_\zs{*}\ge 1$, which may be depend on $n$,
\begin{equation}\label{sec:In.3-00}
\sup_\zs{0\le j\le n}\,\sup_\zs{0\le t\le 1}\vert\phi_\zs{j}(t)\vert\,
\le\,
\phi_\zs{*}
<\infty\,. 
\end{equation}
We extend the functions $\phi_\zs{j}(t)$ by periodicity, i.e.,
 we set $\phi_\zs{j}(t):=\phi_\zs{j}(\{t\})$, where $\{t\}$ is the fractional part of $t\ge 0$.
For example, we can take 
 the trigonometric basis    defined as $\Tr_\zs{1}\equiv 1$ and, for $j\ge 2,$
\begin{equation}\label{sec:In.5}
 \Tr_\zs{j}(x)= \sqrt 2
\left\{
\begin{array}{c}
\cos(2\pi[j/2] x)\, \quad\mbox{for even}\quad j;\\[4mm]
\sin(2\pi[j/2] x)\quad\mbox{for odd}\quad j,
\end{array}
\right.
\end{equation}
where $[x]$ denotes the integer part of $x$.

To estimate the function $S$ we use here the model selection procedure for continuous time regression models from
\cite{KonevPergamenshchikov2012} based on the Fourrier expansion. We recall that  
for any function $S$ from $\L_\zs{2}[0,1]$ we can write 
\begin{equation}\label{sec:In.5+++Fourrier}
S(t)=\sum^{\infty}_\zs{j=1}\,\theta_\zs{j}\,\phi_\zs{j}(t)
\quad\mbox{and}\quad
\theta_\zs{j}= (S,\phi_\zs{j}) = \int_\zs{0}^{1} S(t) \phi_\zs{j}(t)\d t
\,.
\end{equation}
So, to estimate the function $S$ it suffices to estimate the coefficients $\theta_\zs{j}$ and to replace them in this representation by their estimators.
Using the fact that the function $S$ and $\phi_\zs{j}$ are $1$ - periodic we can write that
$$
\theta_\zs{j}=\frac{1}{n} \int_\zs{0}^{n}\, \phi_\zs{j}(t)\,S(t) \d t
\,.
$$
If we replace here the differential $S(t)\d t$ by the stochastic observed differential $\d y_\zs{t}$
we obtain the natural estimate for $\theta_\zs{j}$ on the time interval $[0,n]$ 
\begin{equation}\label{sec:In.7}
\wh{\theta}_\zs{j,n}= \frac{1}{n} \int_\zs{0}^{n}  \phi_\zs{j}(t) \d\,y_\zs{t}\,,
\end{equation}
which
can be represented, in view of the model \eqref{sec:In.1}, as
\begin{equation}\label{sec:In.8}
\wh{\theta}_\zs{j,n}= \theta_\zs{j} + \frac{1}{\sqrt n}\xi_\zs{j,n}\,,
\quad 
\xi_\zs{j,n}= \frac{1}{\sqrt n} I_\zs{n}(\phi_\zs{j})\,.
\end{equation}
Now (see, for example, \cite{IbragimovKhasminskii1981}) we can estimate the function $S$ by the projection estimators, i.e. 
\begin{equation}\label{sec:In.7++pre}
\wh{S}_\zs{m}(t)=\sum^{m}_\zs{j=1}\,\wh{\theta}_\zs{j,n}\,\phi_\zs{j}(t)\,,\quad 0\le t\le 1\,,
\end{equation}
for some number $m\to\infty$ as $n\to\infty$. It should be noted that Pinsker in \cite{Pinsker1981} shows that the projection estimators of the form
\eqref{sec:In.7++pre}
 are not efficient. For obtaining efficient estimation one needs to use weighted least square estimators defined as 
\begin{equation}\label{sec:Mo.1}
\wh{S}_\lambda (t) = \sum_\zs{j=1}^{n} \lambda(j) \wh{\theta}_\zs{j,n} \phi_\zs{j}(t)\,,
\end{equation}
where the coefficients $\lambda=(\lambda(j))_\zs{1\le j\le n}$ belong to some finite set $\Lambda$ from $[0,1]^n$.
 As it is shown in  \cite{Pinsker1981}, in order to obtain efficient estimators, the coefficients $\lambda(j)$ in \eqref{sec:Mo.1} need to be chosen depending on 
the regularity of the unknown function $S$. In this paper we consider the adaptive case, i.e. we assume that the regularity of the function $S$ is unknown.
 In this case we chose the weight coefficients on the basis of  the model selection procedure proposed in \cite{KonevPergamenshchikov2012}
   for the general semi-martingale regression model in continuous time. These coefficients will be obtained later in \eqref{sec:Ga.2}.
 To the end, first we set
 \begin{equation}\label{sec:Mo.2}
\check{\iota}=\#(\Lambda)
\quad\mbox{and}\quad
\vert\Lambda\vert_\zs{*}=1+ \max_\zs{\lambda\in\Lambda}\,\check{L}(\lambda)
\,,
\end{equation}
where $\#(\Lambda)$ is the cardinal number of  $\Lambda$ and $\check{L}(\lambda)=\sum^{n}_\zs{j=1}\lambda(j)$. 
Now, to choose a weight sequence $\lambda$ in the set $\Lambda$ we use the empirical  quadratic risk, defined as
$$
\Er_n(\lambda) = \parallel \wh{S}_\lambda-S\parallel^2,
$$
which in our case is equal to
\begin{equation}\label{sec:Mo.3}
\Er_n(\lambda) = \sum_\zs{j=1}^{n} \lambda^2(j) \wh{\theta}^2_\zs{j,n} -2 \sum_\zs{j=1}^{n} \lambda(j) \wh{\theta}_\zs{j,n}\theta_\zs{j}+ \sum_\zs{j=1}^{\infty} \theta^2_\zs{j}.
\end{equation}
Since the Fourier coefficients $(\theta_\zs{j})_\zs{j\ge\,1}$ are unknown, we replace
the terms $\wh{\theta}_\zs{j,n}\theta_\zs{j,n}$ by  
\begin{equation}\label{sec:Mo.4}
\wt{\theta}_\zs{j,n} = \wh{\theta}^2_\zs{j,n} - \frac{  \wh{\sigma}_\zs{n}}{n}\,,
\end{equation}
where $\wh{\sigma}_\zs{n}$ is an estimate for the variance proxy $\sigma_\zs{Q}$ defined in \eqref{sec:Ex.5}. 
If it is known, we take $\wh{\sigma}_\zs{n}=\sigma_\zs{Q}$; otherwise, we can choose it, for example, as in \cite{KonevPergamenshchikov2012}, i.e.
\begin{equation}\label{sec:Mo.4-1-31-3}
\wh{\sigma}_\zs{n}=
\sum^n_\zs{j=[\sqrt{n}]+1}\,\wh{t}^2_\zs{j,n}\,,
\end{equation}
where $\wh{t}_\zs{j,n}$ are the estimators for the Fourier coefficients with respect to the trigonometric basis \eqref{sec:In.5}, i.e.
\begin{equation}\label{sec:Mo.4-2-31-3}
\wh{t}_\zs{j,n}=\frac{1}{n}
\int^{n}_\zs{0}\,Tr_\zs{j}(t)\d y_\zs{t}\,.
\end{equation}

\noindent 
Finally, in order to choose the weights, we will minimize the following cost function
\begin{equation}\label{sec:Mo.5}
J_n(\lambda)=\sum_\zs{j=1}^{n} \lambda^2(j) \wh{\theta}^2_\zs{j,n} -2 \sum_\zs{j=1}^{n} \lambda(j)\wt{\theta}_\zs{j,n} + \delta\,P_\zs{n}(\lambda),
\end{equation}
where $\delta>0$ is some  threshold which will be specified later and the penalty term is
\begin{equation}\label{sec:Mo.6}
P_\zs{n}(\lambda)= \frac{  \wh{\sigma}_\zs{n} |\lambda|^2}{n}.
\end{equation}
\noindent
We define the model selection procedure as
\begin{equation}\label{sec:Mo.9}
\wh{S}_* = \wh{S}_\zs{\hat \lambda}\,,
\end{equation}
where
\begin{equation}\label{sec:Mo.8}
\wh{\lambda}= \mbox{argmin}_\zs{\lambda\in\Lambda} J_n(\lambda).
\end{equation}
We recall that the set $\Lambda$ is finite so $\hat \lambda$ exists. In the case when $\hat \lambda$ is not unique, we take one of them.

Let us now specify the weight coefficients
$(\lambda(j))_\zs{1\le j\le n}$. Consider, for some fixed $0<\varepsilon<1,$
a numerical grid of the form
\begin{equation}\label{sec:Ga.0}
\cA=\{1,\ldots,k^*\}\times\{\varepsilon,\ldots,m\varepsilon\}\,,
\end{equation}
where $m=[1/\ve^2]$. We assume that both parameters $k^*\ge 1$ and $\varepsilon$ are functions of $n$, i.e.
$k^*=k^*(n)$ and $\ve=\ve(n)$, such that
\begin{equation}\label{sec:Ga.1}
\left\{
\begin{array}{ll}
&\lim_\zs{n\to\infty}\,k^*(n)=+\infty\,,
\quad\lim_\zs{n\to\infty}\,\dfrac{k^*(n)}{\ln n}=0\,,\\[6mm]
&
\lim_\zs{n\to\infty}\,\varepsilon(n)=0
\quad\mbox{and}\quad
\lim_\zs{n\to\infty}\,n^{\check{\delta}}\ve(n)\,=+\infty
\end{array}
\right.
\end{equation}
for any $\check{\delta}>0$. One can take, for example, for $n\ge 2$
\begin{equation}\label{sec:Ga.1-00}
\ve(n)=\frac{1}{ \ln n }
\quad\mbox{and}\quad
k^*(n)=k^{*}_\zs{0}+\sqrt{\ln n}\,,
\end{equation}
where $k^{*}_\zs{0}\ge 0$ is some fixed constant and the threshold $\varsigma^{*}(n)$ is introduced in \eqref{sec:Ex.5}.
 For each $\alpha=(\beta, \l)\in\cA$, we introduce the weight
sequence
$$
\lambda_\zs{\alpha}=(\lambda_\zs{\alpha}(j))_\zs{1\le j\le n}
$$
with the elements
\begin{equation}\label{sec:Ga.2}
\lambda_\zs{\alpha}(j)=\Chi_\zs{\{1\le j<j_\zs{*}\}}+
\left(1-(j/\omega_\alpha)^\beta\right)\,
\Chi_\zs{\{ j_\zs{*}\le j\le \omega_\zs{\alpha}\}},
\end{equation}
where
$j_\zs{*}=1+\left[\ln\upsilon_\zs{n}\right]$, $\omega_\zs{\alpha}=(\d_\zs{\beta}\,\l\upsilon_\zs{n})^{1/(2\beta+1)}$,
$$
\d_\zs{\beta}=\frac{(\beta+1)(2\beta+1)}{\pi^{2\beta}\beta}
\quad\mbox{and}\quad
\upsilon_\zs{n}=n/\varsigma^{*}
\,.
$$
Now we define the set $\Lambda$ 
as
\begin{equation}\label{sec:Ga.3}
\Lambda\,=\,\{\lambda_\zs{\alpha}\,,\,\alpha\in\cA\}\,.
\end{equation}
It will be noted that in this case the cardinal of the set $\Lambda$ is  
\begin{equation}
\label{sec:Ga.1++1--1}
\check{\iota}=k^{*} m\,.
\end{equation}
Moreover,
taking into account that $\d_\zs{\beta}<1$ for $\beta\ge 1$ 
we obtain for the set \eqref{sec:Ga.3}
\begin{equation}
\label{sec:Ga.1++1--2}
 \vert \Lambda\vert_\zs{*}\,
 \le\,1+
\sup_\zs{\alpha\in\cA}
  \omega_\zs{\alpha}
\le 1+(\upsilon_\zs{n}/\ve )^{1/3}\,.
\end{equation}

\begin{remark} \label{Re.sec.Ex.2_2}
Note that the form \eqref{sec:Ga.2}
for the weight coefficients in \eqref{sec:Mo.1}
 was proposed by Pinsker in \cite{Pinsker1981}
 for the efficient estimation in the nonadaptive case, i.e. when the regularity parameters of the function $S$ are known. 
In the adaptive case  these weight coefficients are
  used in \cite{KonevPergamenshchikov2012, KonevPergamenshchikov2015}
   to show the asymptotic efficiency for model selection procedures.
\end{remark}

\bigskip

\bigskip

\section{Main results}\label{sec:Mrs}

In this section we obtain in Theorem \ref{Th.sec:Mrs.1} the non-asymptotic oracle inequality for the quadratic risk \eqref{sec:In.4} for the model selection procedure \eqref{sec:Mo.9} and in 
Theorem \ref{Th.sec:Mrs.2} the non-
asymptotic oracle inequality for the robust risk \eqref{sec:In.6} for the same model selection procedure \eqref{sec:Mo.9}, considered with the coefficients
\eqref{sec:Ga.2}. We give the lower and upper bound for the robust risk in Theorems \ref{Th.sec:Ef.1} and \ref{Th.sec:Ef.2}, and also the optimal convergence rate in  Corollary \ref{Co.sec:Mr.1}.\\

Before stating the non-asymptotic oracle inequality, let us first introduce  the following parameters which will be used for describing the rest term 
in the oracle inequalities. For the renewal density $\rho$ defined in \eqref{sec:Cns.1}
 we set
\begin{equation}\label{sec:Mrs.1-0-1}
\Upsilon(x)=\rho(x)-\frac{1}{\check{\tau}}
\quad\mbox{and}\quad
\Vert\Upsilon\Vert_\zs{1}=\int^{+\infty}_\zs{0}\,\vert\Upsilon(x)\vert\,\d x
\,,
\end{equation}
where 
$\check{\tau}=\E\,\tau_\zs{1}$. In Proposition \ref{Pr.sec:A.1}
we show that $\vert\rho\vert_\zs{*}=\sup_\zs{t\ge 0}\vert\rho(t)\vert<\infty$ and $\Vert\Upsilon\Vert_\zs{1}<\infty$. So, using this, we can introduce the following parameters
\begin{equation}\label{sec:Mrs.1-0+Psi}
\Psi_\zs{Q}=4\varkappa_\zs{Q}\check{\iota}+\left(5+
\frac{4 \check{\iota}}{\sigma_\zs{Q}}\right)\,\left(\sigma_\zs{Q}\,\check{\tau}\,\phi^{2}_\zs{max}\,\Vert\Upsilon\Vert_\zs{1}
+
\phi^{4}_\zs{max}
(1+\sigma^{2}_\zs{Q})^{3}\,\check{\l}
\right)
\end{equation}
and
\begin{equation}\label{sec:Mrs.1-0}
\c^{*}_\zs{Q}=\sigma_\zs{Q}+2\varkappa_\zs{Q}+\sigma_\zs{Q}\,\check{\tau}\,\phi^{2}_\zs{max}\,\Vert\Upsilon\Vert_\zs{1}
+
\phi^{4}_\zs{max}
(1+\sigma^{2}_\zs{Q})^{2}\,\check{\l}
\,,
\end{equation}
where 
$
\check{\l}=
 (4\check{\tau}^{2}+8)\,\Vert\Upsilon\Vert_\zs{1}+5+
 13(1+\check{\tau})^{2}(1+\vert \rho\vert^{2}_\zs{*})(\E Y^{4}_\zs{1})+4\Pi(x^{4})$.
First, let us state the non-asymptotic oracle inequality for the quadratic risk \eqref{sec:In.4} for the model selection procedure \eqref{sec:Mo.9}.
\begin{theorem}\label{Th.sec:OI.1} 
Assume that Conditions $\H_\zs{1})$--$\H_\zs{4})$ hold.
Then, for any $n\ge\,1$ and $0 <\delta< 1/6$,  the estimator of $S$ given in \eqref{sec:Mo.9} satisfies the following oracle inequality
\begin{equation}\label{sec:OI.1}
\mathcal{R}_\zs{Q}(\wh{S}_*,S)\leq\frac{1+3\delta}{1-3\delta} \min_\zs{\lambda\in\Lambda} 
\mathcal{R}_\zs{Q}(\wh{S}_\lambda,S)+
\frac{\Psi_\zs{Q}
 + 10 \vert\Lambda\vert_\zs{*}\, \E_\zs{S} |  \wh{\sigma}_\zs{n} -\sigma_\zs{Q} |}{n\delta}\,.
\end{equation}
\end{theorem}

\bigskip

\noindent Now we  study the estimate \eqref{sec:Mo.4-1-31-3}.
\begin{proposition}\label{Pr.sec:Si.1}
Assume that Conditions $\H_\zs{1})$--$\H_\zs{4})$ hold and that the function  $S(\cdot)$ is continuously
differentiable. 
Then, for any $n\ge 2$,
\begin{equation}\label{sec:Si.3}
\E_\zs{Q,S}|\wh{\sigma}_\zs{n}-\sigma_\zs{Q}|
\le
\frac{
6\|\dot{S}\|^2
+\c^{*}_\zs{Q}}{\sqrt{n}}
\,.
\end{equation}
\end{proposition}

\noindent 
Theorem \ref{Th.sec:OI.1} and Proposition \ref{Pr.sec:Si.1}
implies the following result.

\begin{theorem}\label{Th.sec:Mrs.1}
Assume that Conditions $\H_\zs{1})$--$\H_\zs{4})$ hold
 and that the function $S$ is continuously differentiable. 
  Then, for any $n\ge\, 1 $ and $ 0 <\delta \leq 1/6$, the procedure 
 \eqref{sec:Mo.9}, \eqref{sec:Mo.4-1-31-3}  
satisfies the following oracle inequality
\begin{equation}\label{sec:Mrs.1}
\mathcal{R}_\zs{Q}(\wh{S}_*,S)\leq\frac{1+3\delta}{1-3\delta} \min_\zs{\lambda\in\Lambda} 
\mathcal{R}_\zs{Q}(\wh{S}_\lambda,S)+
\frac{60\wt{\Lambda}_\zs{n}\,
\|\dot{S}\|^2
+\wt{\Psi}_\zs{Q,n}}{n\delta}
\,,
\end{equation}
where $\wt{\Psi}_\zs{Q,n}=10 \wt{\Lambda}_\zs{n}\c^{*}_\zs{Q}+\Psi_\zs{Q}$ and
 $\wt{\Lambda}_\zs{n}=\vert\Lambda\vert_\zs{*}/\sqrt{n}$.
\end{theorem}

\bigskip

\begin{remark} \label{Re.sec.Mo.11_++}
Note that the coefficient $\varkappa_\zs{Q}$ can be estimated as
$\varkappa_\zs{Q}\le (1+\check{\tau}\vert\rho\vert_\zs{*})\sigma_\zs{Q}$.
Therefore,taking into account that $\phi^{4}_\zs{max}\ge 1$, the remainder term in \eqref{sec:Mrs.1}
can be estimated as
\begin{equation}
\label{sec:Mrs.11+ReTerm}
\wt{\Psi}_\zs{Q,n}\le \C_\zs{*}
\left(1+\sigma^{6}_\zs{Q}+\frac{1}{\sigma_\zs{Q}}
\right)(1+\wt{\Lambda}_\zs{n})\check{\iota}\phi^{4}_\zs{max}\,,
\end{equation}
where $\C_\zs{*}>0$ is some constant which is  independent  of the distribution $Q$. 
\end{remark}

Furthermore, let us study the robust risk \eqref{sec:In.6} for the procedure \eqref{sec:Mo.9}. In this case, the distribution family $\cQ_\zs{n}$
consists 
 in all distributions on the Skorokhod space $\cD[0,n]$ of the process
\eqref{sec:Ex.1}  with the parameters satisfying the conditions \eqref{sec:Ex.5} and \eqref{sec:Mrs.5-1}.

Moreover, we assume also that the
number of the weight vectors and
the upper bound for the basis functions  in \eqref{sec:In.3-00} 
may depend on $n\ge 1$, i.e. $\check{\iota}=\check{\iota}(n)$ and $\phi_\zs{*}=\phi_\zs{*}(n)$, such that 
for any $\check{\epsilon}>0$
\begin{equation}\label{sec:Mrs.5-2}
\lim_\zs{n\to\infty}\,\frac{\check{\iota}(n)}{n^{\check{\epsilon}}}=
0
\quad\mbox{and}\quad
\lim_\zs{n\to\infty}\,\frac{\phi_\zs{*}(n)}{n^{\check{\epsilon}}}=
0\,.
\end{equation}

The next result presents the non-asymptotic oracle inequality for the robust risk \eqref{sec:In.6} for the model selection procedure \eqref{sec:Mo.9}, considered with the coefficients
\eqref{sec:Ga.2}. 

\begin{theorem}\label{Th.sec:Mrs.2}
Assume that Conditions $\H_\zs{1})$ -- $\H_\zs{4})$ hold 
and that the unknown function $S$ is continuously differentiable.
Then, for the robust risk defined in  \eqref{sec:In.6} through the distribution family \eqref{sec:Ex.5} -- \eqref{sec:Mrs.5-1},  the procedure \eqref{sec:Mo.9} with the coefficients \eqref{sec:Ga.2}
for any $n\ge\, 1 $ and $ 0 <\delta <1/6$, 
satisfies the following oracle inequality  
\begin{equation}\label{sec:Mrs.6-25.3}
\cR^{*}(\wh{S}_*,S)\leq\frac{1+3\delta}{1-3\delta} \min_\zs{\lambda\in\Lambda}
 \cR^{*}(\wh{S}_\lambda,S)+
\frac{\U^{*}_\zs{n}(S)}{n\delta},
\end{equation}
where the sequence $\U^{*}_\zs{n}(S)>0$ is such that, under the conditions \eqref{sec:Mrs.5-1},
\eqref{sec:Ga.1}  and
 \eqref{sec:Mrs.5-2}, 
for any $r>0$ and $\check{\delta}>0,$
\begin{equation}\label{sec:Mrs.7-25.3}
\lim_\zs{n\to\infty}\,
\sup_\zs{\|\dot{S}\| \le r}
\,
\frac{\U^{*}_\zs{n}(S)}{n^{\check{\delta}}}
=0.
\end{equation}
\end{theorem}

Now we study the asymptotic efficiency for the procedure \eqref{sec:Mo.9} with the coefficients \eqref{sec:Ga.2}, 
  with respect to the robust risk \eqref{sec:In.6} defined by the 
  distribution family \eqref{sec:Ex.5}--\eqref{sec:Mrs.5-1}.  To this end, we assume that the unknown function $S$ in the model
 \eqref{sec:In.1} belongs to the Sobolev ball
\begin{equation}\label{sec:Ef.1}
W^{k}_\zs{r}=\{f\in \,\cC^{k}_\zs{per}[0,1]
\,:\,\sum_\zs{j=0}^k\,\|f^{(j)}\|^2\le \r\}\,,
 \end{equation}
where $\r>0$ 	and $k\ge 1$ are some unknown parameters,
$\cC^{k}_\zs{per}[0,1]$ is the set of
 $k$ times continuously differentiable functions
$f\,:\,[0,1]\to\bbr$ such that $f^{(i)}(0)=f^{(i)}(1)$ for all
$0\le i \le k$. The function class $W^{k}_\zs{r}$ can be written
as an ellipsoid in $\L_\zs{2}[0,1]$, i.e.,
 \begin{equation}\label{sec:Ef.2}
W^{k}_\zs{r}=\{f\in\,\cC^{k}_\zs{per}[0,1]\,:\,
\sum_\zs{j=1}^{\infty}\,a_\zs{j}\,\theta^2_\zs{j}\,\le \r\},
 \end{equation}
where $a_\zs{j}=\sum^k_\zs{i=0}\left(2\pi [j/2]\right)^{2i}$ and $\theta_\zs{j}=\int^{1}_\zs{0}\,f(v)\Tr_\zs{j}(v)\d v$.
We recall that the trigonometric basis 
$(\Tr_\zs{j})_\zs{j\ge 1}$ is defined in
\eqref{sec:In.5}.

Similarly to
\cite{KonevPergamenshchikov2012, KonevPergamenshchikov2015} 
we will  show here that the asymptotic sharp lower bound 
 for the robust risk \eqref{sec:In.6}
is given by
\begin{equation}\label{sec:Ef.3}
\r^{*}_\zs{k}=
\,
\left((2k+1)\r\right)^{1/(2k+1)}\,
\left(
\frac{k}{(k+1)\pi} \right)^{2k/(2k+1)}\,.
\end{equation}

Note that this is  the well-known Pinsker constant
obtained for the nonadaptive filtration problem in ``signal +
small white noise'' model
 (see, for example, \cite{Pinsker1981}).
Let $\Pi_\zs{n}$ be the set of all estimators $\wh{S}_\zs{n}$
measurable with respect to the $\sigma$ - field
$\sigma\{y_\zs{t}\,,\,0\le t\le n\}$
 generated by the process \eqref{sec:In.1}.

The following two results give the lower and upper bound for the robust risk in our case. 

\begin{theorem}\label{Th.sec:Ef.1} Under Conditions \eqref{sec:Ex.5} and \eqref{sec:Mrs.5-1},
 \begin{equation}\label{sec:Ef.4}
\liminf_\zs{n\to\infty}\,
\upsilon
^{2k/(2k+1)}_\zs{n}
 \inf_\zs{\wh{S}_\zs{n}\in\Pi_\zs{n}}\,\,
\sup_\zs{S\in W^{k}_\zs{\r}} \,\cR^{*}_\zs{n}(\wh{S}_\zs{n},S) \ge
\r^{*}_\zs{k}\,,
 \end{equation}
where $\upsilon_\zs{n}=n/\varsigma^{*}$.
\end{theorem}

Note that if the parameters $\r$ and $k$ are known, i.e. for the non adaptive estimation case, then to obtain
the efficient estimation for the "signal+white noise" model Pinsker in \cite{Pinsker1981} proposed to use the estimate 
$\wh{S}_\zs{\lambda_\zs{0}}$ defined in
\eqref{sec:Mo.1} with the weights \eqref{sec:Ga.2} in which
\begin{equation}
\label{sec:Mo.11+l}
\lambda_\zs{0}=\lambda_\zs{\alpha_\zs{0}}
\quad\mbox{and}\quad
\alpha_\zs{0}=(k,\l_\zs{0})\,,
\end{equation}
where $\l_\zs{0}=[\r/\varepsilon ]\varepsilon$. For the model \eqref{sec:In.1} -- \eqref{sec:Ex.1} we show the same result.

\begin{proposition}\label{Th.sec:Ef.33}
The estimator $\wh{S}_{\lambda_\zs{0}}$ satisfies the following asymptotic upper bound
$$
\lim_\zs{n \to \infty } \upsilon^{2k /(2k+1)}_\zs{n}\, \sup_\zs{S\in W^{k}_\zs{\r}} \cR^*_n (\wh{S}_{\lambda_\zs{0}},S) \leq \r^*_\zs{k}\,.
$$
\end{proposition}

\noindent 
For the adaptive estimation we user the model selection procedure \eqref{sec:Mo.9}
with the parameter $\delta$ defined  as a function of $n$ satisfying  
 \begin{equation}\label{sec:Ef.4-01}
\lim_\zs{n}\,\delta_\zs{n}=0
\quad\mbox{and}\quad
\lim_\zs{n}\,n^{\check{\delta}}\,\delta_\zs{n}=0
 \end{equation}
for any $\check{\delta}>0$. For example, we can take $\delta_\zs{n}=(6+\ln n)^{-1}$.

\begin{theorem}\label{Th.sec:Ef.2}
Assume that  Conditions $\H_\zs{1})$--$\H_\zs{4})$ hold true.  
Then  the robust risk defined in  \eqref{sec:In.6} through the distribution family \eqref{sec:Ex.5}--\eqref{sec:Mrs.5-1} 
for
the procedure \eqref{sec:Mo.9} based on the trigonometric basis \eqref{sec:In.5}  with the coefficients
\eqref{sec:Ga.2} and the parameter $\delta=\delta_\zs{n}$ satisfying \eqref{sec:Ef.4-01}
 has
 the following asymptotic upper bound
 \begin{equation}\label{sec:Ef.5}
\limsup_\zs{n\to\infty}\,
\upsilon
^{2k/(2k+1)}_\zs{n}\,
 \sup_\zs{S\in W^{k}_\zs{\r}}\,
\cR^{*}_\zs{n}(\wh{S}_\zs{*},S) \le  
\r^{*}_\zs{k}
\,.
 \end{equation}
\end{theorem}
Theorem~\ref{Th.sec:Ef.1} and Theorem~\ref{Th.sec:Ef.2} allow us to compute the optimal convergence rate.

\begin{corollary}\label{Co.sec:Mr.1}
Under the assumptions of Theorem~\ref{Th.sec:Ef.2}, we have
\begin{equation}\label{sec:Ef.6}
\lim_\zs{n\to\infty}\,
\upsilon
^{2k/(2k+1)}_\zs{n}\,
 \inf_\zs{\wh{S}_\zs{n}\in\Pi_\zs{n}}\,\,
\sup_\zs{S\in W^{k}_\zs{r}} \,\cR^{*}_\zs{n}(\wh{S}_\zs{n},S)
= \r^{*}_\zs{k}\,.
 \end{equation}
\end{corollary}

 \begin{remark}
 \label{Re.sec.Mrs.1}
It is well known that
 the optimal (minimax) risk convergence rate
for the Sobolev ball $W^{k}_\zs{r}$
 is $n^{2k/(2k+1)}$ (see, for example, \cite{Pinsker1981}, \cite{Nussbaum1985}).
We see here that the efficient robust rate is
 $\upsilon^{2k/(2k+1)}_\zs{n}$, i.e., if  the distribution upper bound  $\varsigma^{*}\to 0$  as $n\to\infty,$
 we obtain a faster rate with respect to $n^{2k/(2k+1)}$, and, if $\varsigma^{*}\to \infty$  as $n\to\infty,$
we obtain a slower rate. In the case when $\varsigma^{*}$ is constant, than the robust rate is the same as the classical non robust convergence rate.
\end{remark}

\section{Renewal density}\label{sec:Rtl}

This section is concerned with results related to the renewal measure	
\eqref{sec:Cns.1}. We start with the following  lemma.

\begin{lemma}\label{Le.sec:Rtl.1}
Let $\tau$ be a positive random variable with a density $g$, such that 
$\E e^{\beta \tau} <\infty $ for some $\beta > 0$. Then there exists a constant $\beta_1,$  $0 < \beta_1 \le\beta$ for which,
$$ 
\E e^{(\beta_1+i \omega)\tau} \neq 1 \qquad  \forall \omega \in \bbr\,.
$$
\end{lemma}

\proof We will show this lemma by the contradiction, i.e.
assume  there exists  some sequence of positive numbers going to zero 
$(\gamma_\zs{k} )_\zs{k\ge 1}$ and  a sequence $(w_\zs{k})_\zs{k\ge 1}$ such that  
\begin{equation}\label{sec:Rtl.1-00}
\E e^{(\gamma_k+i \omega_k)\tau} = 1
\end{equation}
for any $k\ge 1$. Firstly assume that 
$\limsup_\zs{k\to\infty}\,w_\zs{k} = +\infty$. Note that in this case, for any $N\geq 1,$
\begin{align*}
\left\vert \int_\zs{0}^{N}\,e^{\gamma_\zs{k}t} \cos(w_\zs{k}t) \,g(t) \d t \right\vert&\le 
\left\vert \int_\zs{0}^{N}\, \cos(w_\zs{k}t) \,g(t) \d t\right\vert\\[2mm]
&+ 
\left\vert\int_\zs{0}^{N}\,(e^{\gamma_\zs{k}t}-1) \cos(w_\zs{k}t) \,g(t) \d t \right\vert\,,
\end{align*}
i.e., in view of Lemma \ref{Le.sec:A.2-00-1}, for any fixed $N\ge 1$ 
$$ 
\limsup_\zs{k\to N}  \int_\zs{0}^{N}\,e^{\gamma_\zs{k}t} \cos(w_\zs{k}t) \,g(t) \d t= 0\,. 
$$
Since  for some $\beta>0$ the integral $ \int_\zs{0}^{+\infty}\,e^{\beta t}\,g(t) \d t<\infty$, we get
$$ 
\lim_\zs{k\to \infty}  \int_\zs{0}^{+\infty}\,e^{\gamma_\zs{k}t} \cos(w_\zs{k}t) \,g(t) \d t= 0\,. 
$$
Let now $ \limsup_\zs{k\to \infty} w_\zs{k}=\omega_\zs{\infty} \neq 0$ and $0<\vert\omega_\zs{\infty}\vert<\infty$. In this case there exists
a sequence 
$(l_k)_\zs{k\ge 1}$ such that $ \lim_\zs{k\to \infty} w_\zs{l_k}=\omega_\zs{\infty}$, i.e.
$$
1=\limsup_\zs{k\to \infty} \E e^{\gamma_\zs{l_k} \tau} \cos(\tau w_\zs{l_k}) =
\E\, \cos(\tau w_\zs{\infty})\,. 
$$
It is clear that, for random variables having density, the last equality is possible if and only if $w_\zs{\infty}=0$.
In this case, i.e. when $\limsup_\zs{k\to \infty} w_\zs{l_k} = 0$, the equation  \eqref{sec:Rtl.1-00} implies
$$
 \limsup_\zs{k\to \infty} \E\, e^{\gamma_\zs{l_k} \tau} \frac{\sin(\tau w_\zs{l_k})}{w_\zs{l_k}} = \E \,\tau =0\,.
 $$
 But, under our conditions, $\E\tau>0$. These contradictions imply the desired result. 
\fdem

\bigskip

\begin{proposition}\label{Pr.sec:A.1}
Let $\tau$ be a positive random variable with the distribution $\eta$
having a density $g$ which satisfies Conditions $\H_\zs{1})$--$\H_\zs{4})$.
  Then the renewal measure 	\eqref{sec:Cns.1} is absolutely continuous with density $\rho$, 
 for which
\begin{equation}\label{sec:A.7-G}
\rho(x)= \frac{1}{\check{\tau}} + \Upsilon(x)\,,
\end{equation}
where $\check{\tau}=\E\tau_\zs{1}$ and $\Upsilon(\cdot)$ is some function defined on  $\bbr_\zs{+}$ with values in $\bbr$   such that
$$
\sup_\zs{x\ge 0}\,x^\gamma\vert\Upsilon(x)\vert <\infty
\quad\mbox{for all}\quad \gamma>0\,.
$$
\end{proposition}
\proof 
First note, that we can represent  the renewal measure 
$\check{\eta}$ as
$\check{\eta}=\eta*\eta_\zs{0}$
and
$\eta_\zs{0}=\sum_\zs{j=0}^{\infty} \eta^{(j)}$. It is clear that in this case the density $\rho$ of  $\check{\eta}$ can be written as 
\begin{equation}\label{sec:Rtl.1}
\rho(x)=\int^{x}_\zs{0}\,
g(x-y)\,
\sum_\zs{n\ge 0}\,g^{(n)}(y)
\d y\,.
\end{equation}
Now we use the arguments proposed in the proof  of  Lemma 9.5 from \cite{Goldie1991}.  
For any $0<\epsilon<1$ we set
\begin{equation}\label{sec:Rtl.01}
\rho_\zs{\epsilon}(x)=\int^{x}_\zs{0}\,
g(x-y)\left(
\sum_\zs{n\ge 0} (1-\epsilon)^{n}\,g^{(n)}(y)-
\frac{(1-\epsilon)}{\check{\tau}}\, g_0(y)
 \right)
 \d y
 -g(x)
 \,,
\end{equation}
where $ g_0(y)= e^{-\epsilon y/\check{\tau}}1_{\{y>0\}}.$ It is easy to deduce that for any $x\in\bbr$
\begin{equation}\label{sec:Rtl.2}
\lim_\zs{\epsilon\to 0}\,\rho_\zs{\epsilon}(x)\,=\,
\rho(x)-\frac{1}{\check{\tau}}\,\int^{x}_\zs{0}\,g(z)\,\d z 
-g(x)\,.
\end{equation}
Moreover,  in view of the condition $\H_\zs{1})$ we obtain that
 the function $\rho_\zs{\epsilon}(x)$ satisfies the condition $\D)$ from 
Section \ref{sec:Four}. So,  through Proposition \ref{Pr.sec:A.2-00}  we get
$$
\rho_\zs{\epsilon}(x+)+\rho_\zs{\epsilon}(x-)=
\frac{1}{\pi}\,\int_\zs{\bbr}\,e^{-ix\theta}\,
\wh{\rho}_\zs{\epsilon}(\theta)\,\d \theta\,,
$$
where $\wh{\rho}_\zs{\epsilon}(\theta)
=\int_\zs{\bbr}\,e^{i\theta x}\rho_\zs{\epsilon}(x)\d x$. 
Note that
$$
\vert\wh{g}(\theta)\vert =  \left
 \vert\int_\zs{\bbr}\,e^{i\theta x}  g(x) \d x  \right
 \vert \leq \int_\zs{\bbr}\,g(x) \d x =1\,,
$$
i.e. for any $0<\epsilon<1$ we have $\vert 1-(1-\epsilon)\wh{g}(\theta)\vert<1$ and therefore
$$ 
\sum_{n=0}^{\infty} (1-\epsilon)^n  (\wh{g}(\theta))^n = \frac{1}{1-(1-\epsilon)\wh{g}(\theta)}\,.
$$
From this and, taking into account that 
\begin{align*}
\wh{g}_0(\theta)= \int_\zs{\bbr}\,e^{i\theta x} g_0(x)\d x= \frac{\check{\tau}}{\epsilon-i\check{\tau}\theta}\,,
\end{align*}
 we obtain 
\begin{align*}
\wh{\rho}_\zs{\epsilon}(\theta)
&= \wh{g}(\theta) \sum_{n=0}^{\infty} (1-\epsilon)^n  (\wh{g}(\theta))^n - \left( \frac{1-\epsilon}{\check{\tau}}\right)
\wh{g}(\theta) \wh{g}_0(\theta)- \wh{g}(\theta)\\[2mm]
&=\wh{g}(\theta)
G_\zs{\epsilon}(\theta)
\quad\mbox{and}\quad
G_\zs{\epsilon}(\theta)
=\frac{1}{1-(1-\epsilon)\wh{g}(\theta)}
-\frac{1-i\check{\tau}\theta}{\epsilon-i\check{\tau}\theta}\,,
\end{align*}
i.e.
\begin{equation}\label{sec:Rtl.3}
\rho_\zs{\epsilon}(x-)
+
\rho_\zs{\epsilon}(x+)
=
\frac{1}{\pi}\,
\int_\zs{\bbr}\,e^{-ix\theta}\,
\wh{g}(\theta)
G_\zs{\epsilon}(\theta)\,\d \theta
\,.
\end{equation}
One can check directly that 
$$
\sup_\zs{0<\epsilon<1,\theta\in\bbr}\vert G_\zs{\epsilon}(\theta)\vert\,<\,\infty\,.
$$
Therefore, using the condition $\H_\zs{3})$ and the
Lebesgue's dominated convergence theorem, we can pass to limit as $\epsilon\to 0$ in \eqref{sec:Rtl.3}, i.e., we obtain that
$$
\rho(x+)
+
\rho(x-)
-\frac{2}{\check{\tau}}\,\int^{x}_\zs{0}\,g(z)\,\d z 
-g(x+)
-g(x-)
=\frac{1}{\pi}\,
\int_\zs{\bbr}\,e^{-ix\theta}\,
\wh{g}(\theta)
G_\zs{0}(\theta)\,\d \theta
\,,
$$
where 
$$
G_\zs{0}(\theta)
=\frac{1}{1-\wh{g}(\theta)}
+
\frac{1-i\check{\tau}\theta}{i\check{\tau}\theta}\,.
$$
Using here again Proposition \ref{Pr.sec:A.2-00} we deduce that

\begin{equation}\label{sec:Rtl.4-00}
\rho(x+)
+
\rho(x-)
=
\frac{2}{\check{\tau}}\,\int^{x}_\zs{0}\,g(z)\,\d z 
+
\frac{1}{\pi}\,
\int_\zs{\bbr}\,e^{-ix\theta}\,
\wh{g}(\theta)
\check{G}(\theta)\,\d \theta
\end{equation}
and
$$
\check{G}(\theta)
=\frac{1}{1-\wh{g}(\theta)}
+
\frac{1}{i\check{\tau}\theta}\,.
$$
Note now that we can represent the density 
\eqref{sec:Rtl.1} as
\begin{align*}
\rho(x)=g*\sum_\zs{n\ge 0}\,g^{(n)}
=\sum_\zs{n\ge 1}\,g^{(n)}(x)
=g(x)+\sum_\zs{n\ge 2}\,g^{(n)}(x)
=:g(x)+\rho_\zs{c}(x)
\end{align*}
and the function $ \rho_\zs{c}(x)$ is continuous for all $x\in\bbr$.  This means that
$$
\wt{\rho}(x)=
\frac{\rho(x+)+\rho(x-)}{2}-\rho(x)
=
\frac{g(x+)+g(x-)}{2}-g(x)
$$
and, therefore,  the condition $\H_\zs{2})$ implies that,
for any $\gamma>0,$
$$
\sup_\zs{x\ge 0}\,x^{\gamma}\,\vert \wt{\rho}(x)\vert\,<\infty.
$$
 Now we can rewrite \eqref{sec:Rtl.4-00}
as
\begin{equation}\label{sec:Rtl.4-01}
\rho(x)
=
\frac{1}{\check{\tau}}\,\int^{x}_\zs{0}\,g(z)\,\d z 
+
\frac{1}{2\pi}\,
\int_\zs{\bbr}\,e^{-ix\theta}\,
\wh{g}(\theta)
\check{G}(\theta)\,\d \theta
-\wt{\rho}(x).
\end{equation}
Taking into account that $\E\,e^{\beta\tau}<\infty$ for some $\beta>0$ we can obtain that
$$
\sup_\zs{x\ge 0}\,x^{\gamma}\,\int^{+\infty}_\zs{x}\,g(z)\,\d z 
<\infty\,.
$$
To study the second term in (\ref{sec:Rtl.4-01}) we will use  Proposition \ref{Pr.sec:A.1-00}. Indeed, the condition  
$\H_\zs{3})$ implies the first limit equality in \eqref{sec:A.7-00}. The second one follows directly from
 Lemma \ref{Le.sec:A.2-00-1}. Therefore, in view of  Proposition \ref{Pr.sec:A.1-00},
 there exists some $\beta^{*}>0$ such that, for any $0\le \beta_\zs{0}\le \beta^{*},$
$$
\int_\zs{\bbr}\,e^{-ix\theta}\,
\wh{g}(\theta)
\check{G}(\theta)\,\d \theta
=
e^{-\beta_\zs{0}x}\,
\int_\zs{\bbr}\,e^{-ix\theta}\,
\wh{g}(\theta-i\beta_\zs{0})
\check{G}(\theta-i\beta_\zs{0})\,\d \theta\,.
$$
Note that, due to Lemma \ref{Le.sec:Rtl.1}, the function $1-\wh{g}(\theta)$ has no zeros on the line 
 $\left\{z\in\bbc\,:\,\mbox{Im}(z)=-\beta_\zs{1} \right\}$. Moreover, one can check
 directly that $\theta=0$ is an isolated zero. So, this means that for any $N>1$
 there can be only finitely many zeros in $\left\{z\in\bbc\,:\,-\beta_\zs{1}<\mbox{Im}(z)<0\,,\,\vert \mbox{Re}(z)\vert<N \right\}$ of the function $1-\wh{g}(\theta).$ Moreover, note that in view of lemma \ref{Le.sec:A.2-00-1} 
for any $r>0$
$$
\lim_\zs{\mbox{Re}(\theta)\to\infty,\vert \mbox{Im}(\theta)\vert\le r}\,\wh{g}(\theta)\,
=0\,.
$$
This means that there exists $N>0$ such that the function $1-\wh{g}(\theta)\neq 0$
for $\theta\in \left\{z\in\bbc\,:\,-\beta_\zs{1}<\mbox{Im}(z)<0\,,\,\vert \mbox{Re}(z)\vert\ge N \right\}$.
So, there can be only finitely many zeros of
the function $1-\wh{g}(\theta)$
 in $\left\{z\in\bbc\,:\,-\beta_\zs{1}<Im(z)<0 \right\}$
for some fixed $0<\beta_\zs{1}<\beta$. Therefore, there exists some $\beta_\zs{0}>0$
for which the function $1-\wh{g}(\theta)$ has no zeros in 
$\left\{z\in\bbc\,:\,-\beta_\zs{0}<Im(z)<0 \right\}$, i.e. the function $\check{G}(\theta)$
will be bounded in this set and we obtain that
$$
\sup_\zs{x\ge 0}\,e^{\beta_\zs{0}x}\,
\left\vert\,
\int_\zs{\bbr}\,e^{-ix\theta}\,
\wh{g}(\theta)
\check{G}(\theta)\,\d \theta
\right\vert
<\infty\,.
$$
This the conclusion follows.  \fdem

Using this proposition we can study the renewal process  $(N_\zs{t})_\zs{t\ge 0}$ introduced in \eqref{sec:Ex.4}.

\begin{corollary}\label{Co.sec:Rtl.1-20-3} Assume that
Conditions $\H_\zs{1})$--$\H_\zs{4})$ hold true.
Then, for any $t>0,$
\begin{equation}\label{sec:Rtl.1-20-3}
\E\,N_\zs{t}\le \vert \rho\vert_\zs{*}\,t
\quad\mbox{and}\quad
\E\,N^{2}_\zs{t}\le\,
\vert \rho\vert_\zs{*}\,t
+
 \,\vert \rho\vert^{2}_\zs{*}\,t^{2}\,.
\end{equation}
\end{corollary}
\proof 
First, by means of Proposition \ref{Pr.sec:A.1}, note that we get 
$$
\E\,N_\zs{t}=\E\,\sum_\zs{k\ge 1}\,\Chi_\zs{\{T_\zs{k}\le t\}}=\int^{t}_\zs{0}\,\rho(v)\,\d v\le \vert\rho\vert_\zs{*}\,t\,.
$$
Regarding the last bound in \eqref{sec:Rtl.1-20-3}, we use the same reasoning as in the previous inequality, i.e., we obtain  
\begin{align*}
\E\,N^{2}_\zs{t}&=\E\,\sum_\zs{k\ge 1}\,\Chi_\zs{\{T_\zs{k}\le t\}}+2\E\,\sum_\zs{k\ge 1}\,\Chi_\zs{\{T_\zs{k}\le t\}}\,\sum_\zs{j=k+1}\,\Chi_\zs{\{T_\zs{j}\le t\}}\\[2mm]
&=\E\,N_\zs{t}+2\E\,\sum_\zs{k\ge 1}\,\Chi_\zs{\{T_\zs{k}\le t\}}\,\Theta(T_\zs{k})=\E\,N_\zs{t}
+\int^{t}_\zs{0}\,\Theta(v)\,\rho(v)\,\d v
\,,
\end{align*}
where, for $0\le v\le t,$ we defined the function $\Theta(v)=\E\,N_\zs{t-v}\le \vert \rho\vert_\zs{*}(t-v)$. \fdem

\section{Stochastic calculus for semi-Markov processes}\label{sec:Smp}

In this section we give some results of stochastic calculus for the process $(\xi_\zs{t})_\zs{t\ge\, 0}$ given in \eqref{sec:Ex.1}, needed all along this paper. As the process $\xi_\zs{t}$ is the combination of a Levy process and a semi-Markov process, these results are not standard and need to be provided.

\begin{lemma}\label{Le.sec:Smp.2} 
Let $f$ and $g$ be any non-random functions from 
$\L_\zs{2}[0,n]$ and $(I_\zs{t}(f))_\zs{t\ge\,0}$ be the process defined in \eqref{sec:In.2}. Then, for any $0 \leq t \leq n$,
\begin{equation}\label{sec:A.06-11-00-00}
\E\,I_\zs{t}(f) I_\zs{t}(g) = \varrho_\zs{1}^2
\,
(f,g)_\zs{t}\,
+ 
\varrho_\zs{2}^2
\,
(f,g\rho)_\zs{t}
\,, 
\end{equation}
where
$(f,g)_\zs{t}=\int_\zs{0}^{t} f(s)\,g(s) \d s$ and
 $\rho$ is the  density defined in \eqref{sec:Cns.1}.
\end{lemma}
\proof
First, note that we can represent the stochastic integral $I_\zs{t}(f)$ as
\begin{equation}
\label{sec:Smp.++1-rpr}
I_\zs{t}(f) =  \varrho_\zs{1} I_\zs{t}^L(f)+ \varrho_\zs{2} I_\zs{t}^z(f)
\,,
\end{equation}
where 
$$
I_\zs{t}^L(f)=\int_\zs{0}^{t}\, f(s) \d L_\zs{s}
\quad\mbox{and}\quad
I_\zs{t}^z(f)=\int_\zs{0}^{t}\, f(s) \d z_\zs{s}
\,.
$$
Note that the mutual covariation
for the martingales $I_\zs{t}^L(f)$ and $I_\zs{t}^L(g)$
(see, for example, \cite{Liptser_Shiryaev_Mart_1986})
may be calculated as
\begin{equation}
\label{sec:Smp.++mutuel_var}
[I^L(f),I^L(g)]_\zs{t}
=\check{\varrho}^2\int^{t}_\zs{0}f(s)g(s)\d s
+(1-\check{\varrho}^2)\sum_\zs{0\le s\le t}f(s)g(s) \left(\Delta \check{L}_\zs{s} \right)^{2}
\,,
\end{equation}
where $\Delta \check{L}_\zs{s}=\check{L}_\zs{s}-\check{L}_\zs{s-}$. Taking into account
 that $\E\,I_\zs{t}^L(f)\,I_\zs{t}^L(g)=\E\,[I^L(f),I^L(g)]_\zs{t}$
and that in view of the first condition in  \eqref{sec:Mcs.2}
 $\Pi(x^{2})=1$, we obtain that
\begin{align}\nonumber
\E\,I_\zs{t}^L(f)\,I_\zs{t}^L(g)
&=
\check{\varrho}^2 \,\int^{t}_\zs{0}\,f(s)g(s)\d s
+ (1-\check{\varrho}^2) \,\Pi(x^{2})\,
\int^{t}_\zs{0}\,f(s)\,g(s)\d s \\ \label{sec:Smp.++mutuel_var.1}
&=
\int^{t}_\zs{0}\,f(s)\,g(s)\d s
\,.
\end{align}
Moreover, note that
\begin{align*}
\E I_\zs{t}^z(f) I_\zs{t}^z(g) & =  \E \left(\sum_\zs{l=1}^{\infty} f(T_\zs{l}) g(T_\zs{l}) Y^2_\zs{l} 1_\zs{\{T_\zs{l} \leq t\}} \right)\\[2mm]
& =  \E \left(\sum_\zs{l=1}^{\infty} f(T_\zs{l}) g(T_\zs{l}) 1_\zs{\{T_\zs{l} \leq t\}} \right) = \int_\zs{0}^{t} f(s) g(s) \rho(s) \d s\,. 
\end{align*}
Hence the conclusion follows.
\fdem

\begin{lemma}\label{Le.sec:Smp.1} 
Assume that Conditions $\H_\zs{1})$--$\H_\zs{4})$ hold true. Then, for any $n\ge 1$ and for any non random function $f$ from $\L_\zs{2}[0,n],$ the stochastic integral \eqref{sec:In.2} exists and satisfies the properties \eqref{sec:In.3} with the coefficient $\varkappa_\zs{Q}$ given in 
\eqref{sec:In.3}.
\end{lemma}
\proof
This lemma follows directly from
Lemma \ref{Le.sec:Smp.2} with $f=g$
and Proposition \ref{Pr.sec:A.1}. 
 \fdem

\begin{lemma}\label{Le.sec:A.06-11-03} 
Let $f$ and $g$ be bounded functions defined on $[0,\infty) \times \bbr.$ Then, for any $k\ge 1,$
$$
\E \left( I_\zs{T_\zs{k-}} (f)\,I_\zs{T_\zs{k-}} (g) \mid \cG \right)= \varrho_\zs{1}^2 
(f\,,\,g)_\zs{T_\zs{k}}+ 
\varrho_\zs{2}^2 \sum_\zs{l=1}^{k-1}\, f(T_\zs{l})\,g(T_\zs{l}),
$$
where  $\cG$ is the $\sigma$-field generated by the sequence $(T_\zs{l})_\zs{l\ge 1}$, i.e., $\cG=\sigma\{T_\zs{l}\,,\,l\ge 1\}$.
\end{lemma}
\proof
Using \eqref{sec:Smp.++1-rpr}, \eqref{sec:Smp.++mutuel_var.1} and, 
taking into account that the process $(L_\zs{t})_\zs{t\ge 0}$ is independent of the $\cG$, we obtain 
$$
\E \left( I_\zs{T_\zs{k-}} (f)\,I_\zs{T_\zs{k-}} (g) \mid \cG \right) 
=\varrho_\zs{1}^2
 (f\,,\,g)_\zs{T_\zs{k}}
+
\E \left( I^{z}_\zs{T_\zs{k-}} (f)\,I^{z}_\zs{T_\zs{k-}} (g) \mid \cG \right) 
\,.
$$
Moreover,
\begin{align*}
\E \left( I^{z}_\zs{T_\zs{k-}} (f)\,I^{z}_\zs{T_\zs{k-}} (g) \mid \cG \right) 
 &=
\E \left( \left(\sum_\zs{l=1}^{k-1}\,f(T_\zs{l})Y_\zs{l}\right) \left(\sum_\zs{l=1}^{k-1}\,g(T_\zs{l})Y_\zs{l}\right)  \mid \cG \right)\\
& =  \sum_\zs{l=1}^{k-1}\, f(T_\zs{l})\,g(T_\zs{l})
 \,.
\end{align*}
This we obtain the desired result. \fdem

\begin{lemma}\label{Le.sec:A.06-11-04}
Assume that Conditions $\H_\zs{1})$--$\H_\zs{4})$ hold true. Then, for any measurable bounded non-random functions $f$ and $g,$ we have
$$
\left\vert
\E\, \int_\zs{0}^{n} I^2_\zs{t-}(f)\, g(t)\, \d m_\zs{t} 
\right\vert
\le 
2
\varrho^{2}_\zs{2}
\vert g\vert_\zs{*}
\vert f\vert^{2}_\zs{*}\,
 \Vert\Upsilon\Vert_\zs{1} 
 \,
n.
$$
\end{lemma}
\proof
Using the definition of the process $(m_\zs{t})_\zs{t\ge 0}$
we can represent this integral as
\begin{align}\nonumber
\int_\zs{0}^{n}\, I^2_\zs{t-}(f)\, g(t)\, \d m_\zs{t}
=&
\sum_\zs{k\ge 1}\,I^2_\zs{T_\zs{k}-}(f)\, g(T_\zs{k})\,Y^{2}_\zs{k}\,
\Chi_\zs{\{T_\zs{k}\le n\}}\\[3mm] \label{sec:A.06-11-00-01}
&-
\int_\zs{0}^{n}\, I^2_\zs{t}(f)\, g(t)\,\rho(t)\, \d t
=: V_\zs{n}-U_\zs{n}
\,.
\end{align}
Note now that
$$
\E\,V_\zs{n}\,
=
\E\,
\sum_\zs{k\ge 1}\,g(T_\zs{k})\,\E\left( I^2_\zs{T_\zs{k}-}(f)\,\mid \cG\right)
\,\Chi_\zs{\{T_\zs{k}\le n\}}\,.
$$
Now, using Lemma \ref{Le.sec:A.06-11-03} we can represent the last expectation as
\begin{equation}\label{sec:A.06-11-01}
\E\,V_\zs{n}\,=
\varrho^{2}_\zs{1}\,\E\,V'_\zs{n}
+
\varrho^{2}_\zs{2}\,\E\,V''_\zs{n},
\end{equation}
where
$$
V^{'}_\zs{n}=
\sum_\zs{k\ge 1}\,g(T_\zs{k})\,
\Vert f\Vert^{2}_\zs{T_\zs{k}}
\,\Chi_\zs{\{T_\zs{k}\le n\}}
\quad\mbox{and}\quad
V^{''}_\zs{n}=
\sum_\zs{k\ge 2}\,g(T_\zs{k})\,
\Chi_\zs{\{T_\zs{k}\le n\}}
\,
\sum^{k-1}_\zs{l=1}
\, f^{2}(T_\zs{l})
\,.
$$
The first term in \eqref{sec:A.06-11-01} can be represented  as
$$
\E\,V^{'}_\zs{n}
=
\int^{n}_\zs{0}\, g(t)\,
\Vert f\Vert^{2}_\zs{t}
\rho(t)\d t\,.
$$
To estimate the last expectation in
\eqref{sec:A.06-11-01}, note that
$$
\E\,
V^{''}_\zs{n}=
\E\,
\sum_\zs{l\ge 1}
\, f^{2}(T_\zs{l})\,
\bar{g}(T_\zs{l})
\Chi_\zs{\{T_\zs{l}\le n\}}
=\int^{n}_\zs{0}\,f^{2}(v)\,\bar{g}(v)\,\rho(v)\d v
\,,
$$
where 
$$
\bar{g}(v)=\E\,\sum_\zs{k\ge 1}\,g(v+T_\zs{k})\,
\Chi_\zs{\{T_\zs{k}\le n-v\}}
=\int^{n}_\zs{v}\,g(t)\,\rho(t-v)\d t
\,.
$$
Moreover, using now the representation \eqref{sec:A.06-11-00-00},
we calculate the expectation of the last term in
\eqref{sec:A.06-11-00-01}
$$
\E\,U_\zs{n}=
\varrho^{2}_\zs{1}
\int_\zs{0}^{n}\, \Vert f\Vert^2_\zs{t}\, g(t)\,\rho(t)\, \d t
+
\varrho^{2}_\zs{2}
\int_\zs{0}^{n}\,
\check{f}(t)
\, g(t)\,\rho(t)\, \d t
\,,
$$
where $\check{f}(t)=\int^{t}_\zs{0}\,f^{2}(s)\,\rho(s)\,\d s$. 
This implies that
$$
\E\, \int_\zs{0}^{n} I^2_\zs{t-}(f)\, g(t)\, \d m_\zs{t} 
=
\varrho^{2}_\zs{2}
 \int_\zs{0}^{n}\, g(t)\,\delta(t) \d t
\,,
$$
where 
$
\delta(t)=
\int_\zs{0}^{t}\,
f^{2}(v)\,
\left( 
\rho(t-v)
-
\rho(t)
\right)\,
\rho(v)\,
 \d v$. Note that, in view of Proposition \ref{Pr.sec:A.1}, the function $\delta$ can be estimated as
 $$
 \vert \delta(t)\vert \le \vert f\vert^{2}_\zs{*}\,
  \vert \rho\vert_\zs{*}\,
 \int_\zs{0}^{t}\,
\left\vert 
\Upsilon(t-v)
-
\Upsilon(t)
\right\vert
\,\d v
\le 
\vert f\vert^{2}_\zs{*}\,
  \vert \rho\vert_\zs{*}\,
  \left( 
  \Vert\Upsilon\Vert_\zs{1}
  +
  t \vert\Upsilon(t) \vert
  \right)\,.
 $$
Therefore,
$$
\left\vert
\E\, \int_\zs{0}^{n} I^2_\zs{t-}(f)\, g(t)\, \d m_\zs{t} 
\right\vert
\le 
2
\varrho^{2}_\zs{2}
\vert g\vert_\zs{*}
\vert f\vert^{2}_\zs{*}\,
 \Vert\Upsilon\Vert_\zs{1} 
 \,
 n
$$
and this finishes the proof. 
\fdem

\begin{lemma}\label{Le.sec:A.05-00} 
Assume that Conditions $\H_\zs{1})$--$\H_\zs{4})$ hold true. 
 Then, for any measurable bounded non-random functions $f$ and $g$, one has
$$ 
\E \int_\zs{0}^{n} I^2_\zs{t-}(f)  I_\zs{t-}(g) g(t) d\xi_\zs{t} = 0.
$$
\end{lemma}
\proof
First, note that
$$
\int_\zs{0}^{n} I^2_\zs{t-}(f)  I_\zs{t-}(g) g(t) d\xi_\zs{t}= \varrho_1  
\int_\zs{0}^{n} I^2_\zs{t}(f)  I_\zs{t}(g) g(t) \d L_\zs{t} 
+ \varrho_\zs{2} \, \int_\zs{0}^{n} I^2_\zs{t-}(f)  I_\zs{t-}(g) g(t) \d z_\zs{t}.
$$
Second, we will show that
\begin{equation}\label{sec:A.05}
\E\, \int_\zs{0}^{n} I^2_\zs{t-}(f)  I_\zs{t-}(g) g(t) \d L_\zs{t}=0\,.
\end{equation}
Using the notations \eqref{sec:Smp.++1-rpr},
 we set
$$
J_\zs{1}=\int_\zs{0}^{n} I^2_\zs{t}(f)  I_\zs{t}^{L}(g) g(t) \d L_\zs{t}
\quad\mbox{and}\quad
J_\zs{2}=\int_\zs{0}^{n} I^2_\zs{t}(f)  I_\zs{t}^{z}(g) g(t) \d L_\zs{t},
$$
we obtain that
\begin{equation}\label{sec:A.05-05-11-00}
\int_\zs{0}^{n}
\, I^2_\zs{t}(f)  I_\zs{t}(g) g(t) \d L_\zs{t}=
\varrho_\zs{1}\,J_\zs{1}
+
\varrho_\zs{2}\,J_\zs{2}\,.
\end{equation}
Now let us  recall  the Novikov  inequalities, \cite{Novikov1975}, also referred to as the Bichteler--Jacod inequalities (see  
\cite{BichtelerJacod1983, MarinelliRockner2014})  providing bound  moments of supremum of purely discontinuous local martingales for
any predictable function  $h$ and any        
 $p\ge 2$ 
\begin{equation}
\label{Novikov++}
\E\sup_{0 \le t\le n} \left\vert\int_{[0,t]\times\bbr} h \;\d(\mu-\nu) \right\vert^{p} \le C^{*}_\zs{p}
 \E\,\check{J}_\zs{p,n}(h)\,,
\end{equation} 
where $C^{*}_\zs{p}$ is some positive constant and
$$
\check{J}_\zs{p,n}(h)=
\left(\int_{[0,n]\times\bbr} h^2 \; \d\nu \right)^{p/2}
 +
 \int_{[0,n]\times\bbr} h^p \, \d\nu 
\,.
$$
By applying this inequality for the non-random function $h=(s,x)=g(s)x$, and, 
recalling that $\Pi(x^{8})<\infty$, we obtain,
$$
\sup_{0 \le t\le n} \E \left\vert I_\zs{t}^{\check{L}}(g) \right\vert^{8} < \infty
\,.
$$
Taking into account that, for any non random square integrated function $f,$ the integral $\left(\int_\zs{0}^{t}\, f(s) \d w_\zs{s}\right)$
is Gaussian with the parameters $\left(0, \int^{t}_\zs{0}\,f^{2}(s)\d s\right)$, we obtain 
$$\sup_{0 \le t\le n} \E \left\vert I_\zs{t}^{L}(g) \right\vert^{8} < \infty.$$
%
Finally, by using the Cauchy inequality, we can estimate for any $0<t\le $ the following expectation as
$$
\E\,(I^{L}_\zs{t}(f))^{4}  (I_\zs{t}^{L}(g))^{2}
<
\sqrt{\E\,(I^{L}_\zs{t}(f))^{8}} \sqrt{\E\,(I^{L}_\zs{t}(f))^{4}}
$$
i.e.,
$$
\sup_\zs{0\le t\le n}\,
\E\,(I^{L}_\zs{t}(f))^{4}  (I_\zs{t}^{L}(g))^{2}
<
\infty\,.
$$
Moreover, taking into account that the 
processes
$(L_\zs{t})_\zs{t\ge 0}$ and $(z_\zs{t})_\zs{t\ge 0}$
are independent,
we obtain that 
$$
\E\,(I^{z}_\zs{t}(f))^{4}  (I_\zs{t}^{L}(g))^{2}
=
\E\,(I^{z}_\zs{t}(f))^{4}  
\E\,(I_\zs{t}^{L}(g))^{2}
=
\Pi(x^{2})  \int^{t}_\zs{0}\,g^{2}(s)\d s
\,
\E\,(I^{z}_\zs{t}(f))^{4}\,.
$$
One can check directly here that, for $t>0,$
$$ 
\E\, |I_\zs{t}^z(f)|^4 
\le \vert f\vert_\zs{*}^{4}\, 
\E\,Y_\zs{1}^4\,
 \E N_\zs{t}^2\,.
 $$
Note that the last bound in Corollary \ref{Co.sec:Rtl.1-20-3}
yields
$
\sup_\zs{0 \le t\le n}\,
\E\, (I_\zs{t}^z(f))^4 
<
\infty$
and, therefore,
$$
\sup_\zs{0 \le t\le n}\,
\E\, (I_\zs{t}(f))^4 (I_\zs{t}^{L}(g))^{2}
<
\infty
\,.
$$
It follows directly that $\E J_\zs{1}=0$. Now we study the last term in \eqref{sec:A.05-05-11-00}. To this end,
first note that similarly to the previous reasoning we obtain that
$$
\E  \int_\zs{0}^{n} (I_\zs{t}^L(f))^2 I_\zs{t}^z(g) g(t) \d L_\zs{t} = 0
\quad\mbox{and}\quad
\E  \int_\zs{0}^{n}  I_\zs{t}^L(f) I_\zs{t}^z(f) I_\zs{t}^z(g) g(t) \d L_\zs{t}=0
\,.
$$
Therefore, to show \eqref{sec:A.05} one needs  to show that
\begin{equation}\label{sec:A.06}
\E \int_\zs{0}^{n} (I_\zs{t}^z(f))^2 I_\zs{t}^z(g) g(t)\,
\d L_\zs{t}= 0\,.
\end{equation}
To check this note that, for any $0<t\le n$
and for any bounded function $f,$
$$
I_\zs{t}^z(f)=\sum^{\infty}_\zs{k=1}\,f(T_\zs{k})\,Y_\zs{k}\,\Chi_\zs{\{T_\zs{k}\le t\}}
=
\sum^{N_\zs{n}}_\zs{k=1}\,f(T_\zs{k})\,Y_\zs{k}\,\Chi_\zs{\{T_\zs{k}\le t\}}\,,
$$
i.e.,
$$
\int_\zs{0}^{n} (I_\zs{t}^z(f))^2 I_\zs{t}^z(g) g(t)\,
\d L_\zs{t}
=
\sum^{N_\zs{n}}_\zs{k=1}\,
\sum^{N_\zs{n}}_\zs{l=1}\,
\sum^{N_\zs{n}}_\zs{j=1}\,
f(T_\zs{k})\,
f(T_\zs{l})\,
g(T_\zs{j})\,Y_\zs{j}
Y_\zs{l}\,Y_\zs{k}\,
I_\zs{klj}
\,,
$$
where
$$
I_\zs{klj}=
\int_\zs{0}^{n} 
\Chi_\zs{\{T_\zs{k}\le t\}}
\Chi_\zs{\{T_\zs{l}\le t\}}
\Chi_\zs{\{T_\zs{j}\le t\}}
\d L_\zs{t}\,.
$$
Taking into account that the  $(L_\zs{t})_\zs{t\ge 0}$ is independent of the field $\cG_\zs{z}=\sigma\{z_\zs{t}\,,t\ge 0\},$
we obtain that $\E\left( I_\zs{klj}\vert \cG_\zs{z} \right)=0$. Therefore,
\begin{eqnarray*}
&&\E \int_\zs{0}^{n} (I_\zs{t}^z(f))^2 I_\zs{t}^z(g) g(t)\, \d L_\zs{t} \\
=&& \E 
\sum^{N_\zs{n}}_\zs{k=1}\,
\sum^{N_\zs{n}}_\zs{l=1}\,
\sum^{N_\zs{n}}_\zs{j=1}\,
f(T_\zs{k})\,
f(T_\zs{l})\,
g(T_\zs{j})\,Y_\zs{j}
Y_\zs{l}\,Y_\zs{k}\,
\E\left( I_\zs{klj}\vert \cG_\zs{z} \right)
=0.
\end{eqnarray*}
So, we obtain  \eqref{sec:A.06} and hence the proof is achieved. \fdem

\section{Properties of the regression model \eqref{sec:In.1}}\label{sec:Prsm}


In order to prove the oracle inequalities we  need to study  the conditions
introduced in \cite{KonevPergamenshchikov2012} for the general semi-martingale model \eqref{sec:In.1}. 
To this end we set for any $x\in\bbr^{n}$ the functions
\begin{equation}\label{sec:Prsm.1}
B_\zs{1,Q,n}(x)=
 \sum_\zs{j=1}^{n} x_\zs{j} \, 
 \left( \E_\zs{Q}\xi^2_\zs{j,n} - \sigma_\zs{Q}\right)
 \quad\mbox{and}\quad
B_\zs{2,Q,n}(x)=
 \sum_\zs{j=1}^{n}\,x_\zs{j}\,\wt{\xi}_\zs{j,n}
\,,
\end{equation}
where $\sigma_\zs{Q}$ is defined in
\eqref{sec:Ex.5}
and $\wt{\xi}_\zs{j,n}=\xi^2_\zs{j,n}- \E_\zs{Q}\xi^2_\zs{j,n}$.

\begin{proposition}\label{Pr.sec:OI.2} 
Assume that Conditions $\H_\zs{1})$--$\H_\zs{4})$ hold. Then 
\begin{equation}\label{sec:OI.1-14.3.1}
\sup_\zs{x\in [-1,1]^{n}}\,\vert B_\zs{1,Q,n}(x)\vert
\le 
\C_\zs{1,Q,n}
\,,
\end{equation}
where $\C_\zs{1,Q,n}=
\sigma_\zs{Q}\,\check{\tau}\,\phi^{2}_\zs{max}\,\Vert\Upsilon\Vert_\zs{1}$.
\end{proposition}
\proof
First, note that from
\eqref{sec:Smp.++1-rpr}
we have
$$
\xi_\zs{j,n} = \frac{\varrho_\zs{1}}{\sqrt n} I_n^L (\phi_\zs{j})  + \frac{\varrho_\zs{2}}{\sqrt n} I_n^z (\phi_\zs{j})\,.
$$
So, using \eqref{sec:Smp.++mutuel_var.1}
we can write that
\begin{equation}\label{sec:Pr.19}
\E{\xi^2_\zs{j,n}} = \frac{\varrho^2_1}{ n} \int_\zs{0}^{n} \phi^2_\zs{j}(t) \d\,t + 
\frac{\varrho^2_2}{ n} \E \sum_\zs{l=1}^\infty \phi^2_\zs{j}(T_\zs{l}) 1_\zs{\{T_\zs{l} \leq n\}}\,.
\end{equation}
Proposition \ref{Pr.sec:A.1} implies
\begin{align*}
\E \sum_\zs{l=1}^\infty \phi^2_\zs{j}(T_\zs{l}) 1_\zs{\{T_\zs{l} \leq n\}}\,&=\,
\int^{n}_\zs{0}\,\phi^2_\zs{j}(x)\,\rho(x)\d x\\
&=\frac{1}{\check{\tau}}\,\int_\zs{0}^{n} \phi^2_\zs{j}(x) \d\,x\,
+\int_\zs{0}^{n} \phi^2_\zs{j}(x)\Upsilon(x) \d\,x\,.
\end{align*}
Note that  $\int_\zs{0}^{n} \phi^2_\zs{j}(t) \d\,t=n$. So, in view of the condition
\eqref{sec:In.3-00},
 we  obtain 
\begin{equation}\label{sec:Pr.19-01}
\left\vert 
\E{\xi^2_\zs{j,n}}
-
\sigma_\zs{Q} 
\right\vert
=\frac{\varrho^2_2}{n}\,
\left\vert \int_\zs{0}^{n} \phi^2_\zs{j}(x)\Upsilon(x) \d\,x
\right\vert
\le\, 
\frac{\varrho^2_2}{n}\,
\phi^2_\zs{max}\,
\Vert\Upsilon\Vert_\zs{1}\,.
\end{equation}
Estimating here $\varrho^2_2$ by $\sigma_\zs{Q}\check{\tau}$ we obtain the inequality  \eqref{sec:OI.1-14.3.1} and hence the conclusion follows. 
\fdem

\begin{proposition}\label{Pr.sec:OI.3} 
Assume that Conditions $\H_\zs{1})$--$\H_\zs{4})$ hold.  Then
\begin{equation}\label{sec:OI.1-14.3.1-2}
\sup_\zs{\vert x\vert\le 1}\,
\E_\zs{Q}\,
B^{2}_\zs{2,Q,n}(x)
\,
\le\,\C_\zs{2,Q,n}, 
\end{equation}
where
$\C_\zs{2,Q,n}=\phi^{4}_\zs{max}
(1+\sigma^{2}_\zs{Q})^{3}\,\check{\l}
$
and $\check{\l}$ is given in \eqref{sec:Mrs.1-0}.
\end{proposition}
 \proof
By Ito's formula  one gets
\begin{equation}
\label{sec:Pr.4-11-00}
\d I^2_\zs{t}(f)= 2 I_\zs{t-}(f) \d I_\zs{t}(f)+ \varrho_\zs{1}^2\check{\varrho}^{2}\, f^2(t) \d\,t + \sum_\zs{0\leq s \leq t}f^2(s) (\Delta \xi^{d}_\zs{s})^2\,,
\end{equation}
where $\xi^{d}_\zs{t}=\varrho_\zs{3}\,\check{L}_\zs{t}+\varrho_\zs{2}z_\zs{t}$ and 
$\varrho_\zs{3}=\varrho_\zs{1}\sqrt{1-\check{\varrho}^{2}}$.  Taking into account that the processes 
$(\check{L}_\zs{t})_\zs{t\ge 0}$
and $(z_\zs{t})_\zs{t\ge 0}$ are independent and the time of jumps $T_\zs{k}$ defined in \eqref{sec:Ex.4}
has a density,
we have $\Delta z_\zs{s} \Delta \check{L}_\zs{s}=0$ a.s. for any $s\ge 0$.
Therefore, we can rewrite the differential 
\eqref{sec:Pr.4-11-00} as
\begin{align}\nonumber
\d I^2_\zs{t}(f)=& 2 I_\zs{t-}(f) \d I_\zs{t}(f)+ \varrho_\zs{1}^2\check{\varrho}^{2}\, f^2(t) \d\,t 
\\[2mm] \label{sec:Pr.4-11-00++1}
&+ 
\varrho^{2}_\zs{3}\d \sum_\zs{0\leq s \leq t}f^2(s) (\Delta \check{L}_\zs{s})^2
+
\varrho^{2}_\zs{2}\d \sum_\zs{0\leq s \leq t}f^2(s) (\Delta z_\zs{s})^2
\,.
\end{align}
From Lemma \ref{Le.sec:Smp.1}  it follows that
$$
\E I^2_\zs{t}(f) = \varrho_\zs{1}^{2}\,
\int_0^t f^2(s) ds + \varrho_\zs{2}^2 \int_0^t f^2(s) \rho(s)  \d s
\,.
$$
 Therefore, putting
\begin{equation}
\label{sec:Pr.4-11-01}
\wt{I}_\zs{t}(f) = I_\zs{t}^2(f) -\E I^2_\zs{t}(f)\,,
\end{equation}
 we obtain
$$
\d\wt{I}_\zs{t}(f)= 2 I_\zs{t-}(f) f(t) d\xi_\zs{t}
+
f^2(t) \d \wt{m}_\zs{t}\,,\quad 
\wt{m}_\zs{t}=\varrho^{2}_\zs{3}\check{m}_\zs{t}+\varrho^{2}_\zs{2}m_\zs{t}\,,
$$
where 
$
\check{m}_\zs{t}= \sum_\zs{0\leq s \leq t}(\Delta \check{L}_\zs{s})^2 - t$ and
$
m_\zs{t} = \sum_\zs{0\leq s \leq t}(\Delta z_\zs{s})^2 -
 \int_\zs{0}^{t} \rho(s) \d s$.
For any non-random vector $x= (x_\zs{j})_\zs{1\le j\le n}$ with $\sum^{n}_\zs{j= 1} x^{2}_\zs{j} \leq 1$, we set
\begin{equation}
\label{sec:Pr.4-11-02}
\bar{I}_\zs{t}(x)= \sum_\zs{j=1}^{n} x_\zs{j} \wt{I}_\zs{t}(\phi_\zs{j}).
\end{equation}
Denoting 
\begin{equation}
\label{sec:Pr.4-11-02++}
A_\zs{t} (x) = \sum^{n}_\zs{j=1} x_\zs{j} I_\zs{t}(\phi_\zs{j}) \phi_\zs{j}(t)
\quad\mbox{and}\quad
B_\zs{t}(x) = \sum^{n}_\zs{j=1} x_\zs{j} \phi^2_\zs{j}(t)\,,
\end{equation}
we get the following stochastic differential equation for \eqref{sec:Pr.4-11-02}
$$ 
\d\bar{I}_\zs{t}(x) = 2 A_\zs{t-} (x) \d\xi_\zs{t}
+
 B_\zs{t}(x) \d \wt{m}_\zs{t} 
\,, \quad 
\bar{I}_\zs{0}(x)=0\,.
$$
Applying the Ito's formula one obtains
\begin{align}
\nonumber
\E\,\bar{I}_n^2(x)=& 2\E \int_\zs{0}^{n} \bar{I}_\zs{t-}(x) \d\bar{I}_\zs{t}(x) + 4 \varrho_\zs{1}^2\check{\varrho}^{2} 
\E\,\int_\zs{0}^{n}  A_\zs{t}^2(x) \d\,t\\[2mm] \label{sec:Pr.4-11-03}
 &+\varrho^{2}_\zs{3}\
\E\, \check{D}_\zs{n}(x)
 + 
\varrho_\zs{2}^2\,\E\,D_\zs{n}(x)\,,
\end{align}
where $\check{D}_\zs{n}(x)=\sum_\zs{0\le t\le n}\left(2A_\zs{t-}(x)\Delta \check{L}_\zs{t}+\varrho^{2}_\zs{3} B_\zs{t}(x)(\Delta \check{L}_\zs{t})^{2}\right)^{2}$ and\\
$
D_\zs{n}(x)=
\sum_\zs{k=1}^{+\infty}
\left(
 2A_\zs{T_\zs{k}-}(x) Y_\zs{k}
 +
\varrho_\zs{2}  B_\zs{T_\zs{k}-}(x) Y_\zs{k}^2 
\right)^{2}
\,\Chi_\zs{\{T_\zs{k} \leq n\}}\,.
$
Let us now show that 
\begin{equation}\label{sec:Pr.31}
\left\vert
\E\,
\int_\zs{0}^{n} \bar{I}_\zs{t-}(x) d\bar{I}_\zs{t}(x) 
\right\vert
\le 
2\,\varrho^{4}_\zs{2}\phi^{3}_\zs{max}\,\Vert\Upsilon\Vert_\zs{1}\,
 n^{2}\,.
\end{equation}
To this end, note that
\begin{align*}
\int_\zs{0}^{n} 
\bar{I}_\zs{t-}(x) \d\bar{I}_\zs{t}(x) =& 
2 \sum_\zs{1\le j,l\le\,n} x_\zs{j} x_\zs{l} 
\int_\zs{0}^{n} \wt{I}_\zs{t-}(\phi_\zs{j}) \,I_\zs{t-}(\phi_\zs{l})\phi_\zs{l}(t) \d \xi_\zs{t}
\\[2mm] 
&
+
 \sum_\zs{j=1}^{n} x_\zs{j} \int_\zs{0}^{n} \wt{I}_\zs{t-}(\phi_\zs{j}) B_\zs{t}(x) 
\d \wt{m}_\zs{t} 
\,.
\end{align*}
Using here  Lemma \ref{Le.sec:A.05-00}, we get
$
\E\,
\int_\zs{0}^{n}\,
 \wt{I}_\zs{t-}(\phi_\zs{j})\,
 I_\zs{t-}(\phi_\zs{i}) \phi_\zs{i}(t)\d\xi_\zs{t}
=0$.
Moreover, the process $(\check{m}_\zs{t})_\zs{t\ge 0}$ is a martingale, i.e. $\E\, \int_\zs{0}^{n} \wt{I}_\zs{t-}(\phi_\zs{j}) B_\zs{t}(x) 
\d \wt{m}_\zs{t}=0$. Therefore,
$$
\E\,
\int_\zs{0}^{n} \bar{I}_\zs{t-}(x) d\bar{I}_\zs{t}(x) 
=\varrho_\zs{2}^2 
\sum_\zs{j=1}^{n} x_\zs{j}\E\, \int_\zs{0}^{n} \wt{I}_\zs{t-}(\phi_\zs{j}) B_\zs{t}(x) 
\d m_\zs{t}\,. 
$$
Taking into account here that
 for any non-random bounded
function $f$ 
$$
\E\,\int_\zs{0}^{n} f(t) \d m_\zs{t}=0,
$$
we obtain 
$\E\int_\zs{0}^{n} \wt{I}_\zs{t-}(\phi_\zs{j}) \,B_\zs{t}(x)\, \d m_\zs{t}
=
\E\,\int_\zs{0}^{n} I^2_\zs{t-} (\phi_\zs{j})\,B_\zs{t}(x)\, \d m_\zs{t}$.
So,
Lemma \ref{Le.sec:A.06-11-04}
yields
\begin{align*}
\left\vert
\E\int_\zs{0}^{n} \wt{I}_\zs{t-}(\phi_\zs{j}) \,B_\zs{t}(x) \d m_\zs{t}\right\vert & =
\left\vert
\sum_\zs{l=1}^{n} x_\zs{l} \E\,\int_\zs{0}^{n} I^2_\zs{t-} (\phi_\zs{j}) \phi^2_\zs{l}(t) \d m_\zs{t}
\right\vert \\[2mm]
&\le \,2\,\varrho^{2}_\zs{2}\phi^{3}_\zs{max}\,\Vert\Upsilon\Vert_\zs{1}\,
\sum_\zs{l=1}^{n} \vert x_\zs{l}\vert n\,.
\end{align*}
Therefore,
\begin{align*}
\left\vert
\E\,
\int_\zs{0}^{n} \bar{I}_\zs{t-}(x) d\bar{I}_\zs{t}(x) 
\right\vert
&\le 
2\,\varrho^{4}_\zs{2}\phi^{3}_\zs{max}\,\Vert\Upsilon\Vert_\zs{1}\,
 n
\sum_\zs{1\le l,j\le\,n}\, \vert x_\zs{l} \vert\, 
\vert x_\zs{j} \vert\\[2mm]
&=
2\,\varrho^{4}_\zs{2}\phi^{3}_\zs{max}\,\Vert\Upsilon\Vert_\zs{1}\, n
\left(
\sum_\zs{l=1}^{n}\, \vert x_\zs{l} \vert\right)^{2}.
\end{align*}
Taking into account here that $
\left(
\sum^{n}_\zs{l=1}\, \vert x_\zs{l} \vert\right)^{2}
\le\,
n
\sum_\zs{l\ge\,1}\, x^{2}_\zs{l}\le n$, we obtain \eqref{sec:Pr.31}.  
Reminding that 
$\Pi(x^{2})=1$ we can calculate
 directly that
\begin{equation}\label{sec:Pr.31--}
\E\,\check{D}_\zs{n}(x)
=4\,\E\int^{n}_\zs{0}\,A^{2}_\zs{t}(x)\d t
+\varrho^{4}_\zs{3}\,\Pi(x^{4})
\int^{n}_\zs{0}\,B^{2}_\zs{t}(x)\d t
 \,.
\end{equation}
Note that, thanks to Lemma \ref{Le.sec:Smp.2}, 
we obtain that 
\begin{align*}
\E \int_\zs{0}^{n} A_\zs{t}^2(x) \d\,t 
& = 
\sum_\zs{i,j} x_\zs{i} x_\zs{j} \int_\zs{0}^{n}  \phi_\zs{i}(t)\phi_\zs{j}(t)  \E I_\zs{t}\phi_\zs{i}(t)I_\zs{t}\phi_\zs{j}(t) \d\,t \\[2mm]
& =
\sum_\zs{i,j} x_\zs{i} x_\zs{j} \int_\zs{0}^{n}  \phi_\zs{i}(t)\phi_\zs{j}(t)
\,\int_\zs{0}^{t}\, \phi_\zs{i}(v)\phi_\zs{j}(v)(\varrho^{2}_\zs{1}+\varrho^{2}_\zs{2}\rho(v)) \d v\\[2mm]
& = \frac{1}{2}\varrho^{2}_\zs{1}\,
\sum_\zs{i,j} x_\zs{i} x_\zs{j}\left( \int_\zs{0}^{n}  \phi_\zs{i}(t)\phi_\zs{j}(t)\d t\right)^{2}
+\varrho^{2}_\zs{2}\,A_\zs{1,n}(x)\\
& =
\frac{n^{2}}{2}\varrho^{2}_\zs{1}
+
\varrho^{2}_\zs{2}\,A_\zs{1,n}(x)\,,
\end{align*}
where $A_\zs{1,n}(x)=\sum_\zs{i,j} x_\zs{i} x_\zs{j} \int_\zs{0}^{n}  \phi_\zs{i}(t)\phi_\zs{j}(t)
\,\left(\int_\zs{0}^{t}\, \phi_\zs{i}(v)\phi_\zs{j}(v)\,\rho(v) \d v\right)\d t$.
This term can be estimated 
through  Proposition \ref{Pr.sec:A.1} as
\begin{align*}
\left\vert A_\zs{1,n}(x)\right\vert&=
\left\vert \frac{n^{2}}{2\check{\tau}}\,+\,
\sum_\zs{i,j} x_\zs{i} x_\zs{j} 
\int_\zs{0}^{n}  \phi_\zs{i}(t)\phi_\zs{j}(t)
\,\left(\int_\zs{0}^{t}\, \phi_\zs{i}(v)\phi_\zs{j}(v)\,\Upsilon(v) \d v\right)\d t
\right\vert
\\[2mm]
&\le 
 \frac{n^{2}}{2\check{\tau}}\,+n \,\phi^{4}_\zs{max}\,
 \Vert\Upsilon\Vert_\zs{1}\,
\sum_\zs{i,j} \vert x_\zs{i}\vert \vert x_\zs{j}\vert 
\le\,
 \left( 
 \frac{1}{2\check{\tau}}\,+\,\phi^{4}_\zs{max}\,
 \Vert\Upsilon\Vert_\zs{1}
\right)
 n^{2}\,.
\end{align*}
So, reminding that $\sigma_\zs{Q}=\varrho^{2}_\zs{1}+\varrho^{2}_\zs{2}/\check{\tau}$ and that $\phi_\zs{max}\ge 1$,
 we obtain that 
\begin{align}\nonumber
\E \int_\zs{0}^{n} A_\zs{t}^2(x) \d\,t 
&\le 
\left(
\frac{\sigma_\zs{Q}}{2}
+\,\phi^{4}_\zs{max}\,
 \Vert\Upsilon\Vert_\zs{1}
\right)
 n^{2}\\ \label{sec:Pr.33}
 &\le 
\left(
\frac{1}{4}
+
 \Vert\Upsilon\Vert_\zs{1}
\right)
\phi^{4}_\zs{max}\,
(1+\sigma^{2}_\zs{Q})
\,
 n^{2}
 \,.
\end{align}
Taking into account that
\begin{equation}\label{sec:Pr.33++}
\sup_\zs{t\ge 0}
B^2_\zs{t}(x)\,\le\,\phi^{4}_\zs{max}\,
\left( 
\sum^{n}_\zs{j=1} \vert x_\zs{j}\vert
\right)^{2}
\le 
\,\phi^{4}_\zs{max}\,
n\,,
\end{equation}
that $\phi_\zs{max}\ge 1,$
and that $\varrho^{4}_\zs{1}\le \sigma^{2}_\zs{Q}$
we estimate the expectation in \eqref{sec:Pr.31--}
as
\begin{equation}
\label{sec:Pr.31--++0}
\E\check{D}_\zs{n}\le 4 \phi^{4}_\zs{max}(1+\sigma^{2}_\zs{Q})
\left(1+\Vert\Upsilon\Vert_\zs{1}+\Pi(x^{4}) \right)\,n^{2}
\,.
\end{equation}

Moreover,
 taking into account that the random variable $Y_\zs{k}$ is independent of  $A_\zs{T_\zs{k^-}}(x)$ and of the field 
 $\cG=\sigma\{T_\zs{j}\,,\,j\ge 1\}$ and that $\E\,\left(A_\zs{T_\zs{k^-}}(x)\,\vert \cG\right) =0$,
 we get
\begin{align*}
\E& \sum_\zs{k=1}^{+\infty} B_\zs{T_\zs{k}-}(x)\,A_\zs{T_\zs{k^-}}(x) Y_\zs{k}^3 1_\zs{\{T_\zs{k} \leq n\}} =  
\sum_\zs{k=1}^{+\infty}\E\E \left( B_\zs{T_\zs{k}-}(x)\,A_\zs{T_\zs{k^-}}(x) Y_\zs{k}^3 1_\zs{\{T_\zs{k} \leq n\}}\vert \cG \right) \\
& = \E\,Y_\zs{1}^3\,\E\sum_\zs{k=1}^{+\infty} B_\zs{T_\zs{k}-}(x)  1_\zs{\{T_\zs{k} \leq n\}}\,  \E (A_\zs{T_\zs{k^-}}(x)\vert \cG) =0\,.
\end{align*}
 Therefore, 
\begin{equation}\label{sec:Pr.32}
\E \,D_n(x) =  \varrho^{2}_\zs{2} \E Y_1^4 D_\zs{1,n}(x) + 4 D_\zs{2,n}(x),
\end{equation}
where
$$ 
D_\zs{1,n}(x) =
 \sum_\zs{k=1}^{+\infty} \E\,
B^2_\zs{T_\zs{k}-}(x)1_\zs{\{T_\zs{k} \leq n\}}
\quad\mbox{and}\quad
D_\zs{2,n}(x)= \sum_\zs{k=1}^{+\infty} \E\,A^2_\zs{T_\zs{k^-}}(x)\Chi_\zs{\{T_\zs{k} \leq n\}}\,.
$$
Using the bound \eqref{sec:Pr.33++}
we can estimate 
the term  $D_\zs{1,n}$ as
$
D_\zs{1,n}(x)  \le \phi^{4}_\zs{max} n \E\,N_\zs{n}$.
Using here Corollary  \ref{Co.sec:Rtl.1-20-3}, we obtain
\begin{equation}\label{sec:Pr.34}
D_\zs{1,n}(x)  \le\,\vert\rho\vert_\zs{*}
\phi^{4}_\zs{max} 
\,n^{2}\,.
\end{equation}
Now, to estimate the last term in
\eqref{sec:Pr.32}, note that the process $A_\zs{t}(x)$ can be rewritten as
\begin{equation}\label{sec:Pr.35}
A_\zs{t}(x) = \int_\zs{0}^{t} Q_x(t,s) d\xi_\zs{s},
\quad\mbox{with}\quad 
Q_x(t,s) = \sum^{n}_\zs{j=1} x_\zs{j} \phi_\zs{j}(s)  \phi_\zs{j}(t).
\end{equation}
Applying Lemma \ref{Le.sec:A.06-11-03}
again, we obtain for any $k\ge 1$
$$
\E\,\left( 
A^2_\zs{T_\zs{k^{-}}}(x)\vert \cG
\right)
=
\varrho^{2}_\zs{1}\,
\int_\zs{0}^{T_\zs{k}}\,Q^{2}_\zs{x}(T_\zs{k},s) \d s
+
\varrho^{2}_\zs{2}\,
\sum^{k-1}_\zs{j=1}\,
Q^{2}_\zs{x}(T_\zs{k},T_\zs{j})
\,.
$$
So, we can represent the last term in \eqref{sec:Pr.32} as 
\begin{equation}\label{sec:Pr.35-00-1}
D_\zs{2,n}=
\varrho^{2}_\zs{1}\,D^{(1)}_\zs{2,n}
+
\varrho^{2}_\zs{2}\,D^{(2)}_\zs{2,n}
\,,
\end{equation}
where
$$
D^{(1)}_\zs{2,n}=\sum_\zs{k=1}^{+\infty} \E\,
\Chi_\zs{\{T_\zs{k} \leq n\}}\,
\int_\zs{0}^{T_\zs{k}}\,Q^{2}_\zs{x}(T_\zs{k},s) \d s
$$
 and
$$
D^{(2)}_\zs{2,n}=\sum_\zs{k=1}^{+\infty} \E\,
\Chi_\zs{\{T_\zs{k} \leq n\}}\,
\sum^{k-1}_\zs{j=1}\,
Q^{2}_\zs{x}(T_\zs{k},T_\zs{j})\,.
$$
Thanks to Proposition \ref{Pr.sec:A.1}
we obtain
$$
D^{(1)}_\zs{2,n}=
\int^{n}_\zs{0}\,\left( 
\int_\zs{0}^{t}\,Q^{2}_\zs{x}(t,s) \d s
\right)\,
\rho(t)
\,\d t
\,\le\,
\vert\rho\vert_\zs{*}\,
\int^{n}_\zs{0}\,
\int_\zs{0}^{n}\,Q^{2}_\zs{x}(t,s) \d s
\,\d t
\,.
$$
In view of  the definition of $Q_\zs{x}$ in \eqref{sec:Pr.35},
we can rewrite the last integral as
\begin{align*}
\int_\zs{0}^{n}\,Q^{2}_\zs{x}(t,s) \d s
&=
\sum_\zs{1\le i,j\le n}\,x_\zs{i}\,x_\zs{j}\,
\phi_\zs{i}(t)\,\phi_\zs{j}(t)\,
\int^{n}_\zs{0}\,
\phi_\zs{i}(s)\,\phi_\zs{j}(s)\,\d s\\[2mm]
&
=
\sum^{n}_\zs{i=1}\,x^{2}_\zs{i}\,
\phi^{2}_\zs{i}(t)\,
\int^{n}_\zs{0}\,
\phi^{2}_\zs{i}(s)\,\d s
=\,n\,
\sum^{n}_\zs{i=1}\,x^{2}_\zs{i}\,
\phi^{2}_\zs{i}(t)
\,.
\end{align*}
Since $\sum^{n}_\zs{j=1}\,x^{2}_\zs{j}\le 1$, we obtain that,
\begin{equation}\label{sec:Pr.35-00-1-2}
\int_\zs{0}^{n}\,Q^{2}_\zs{x}(t,s) \d s\,
\le\,
 \phi^{2}_\zs{max}\,
n
\quad\mbox{and}\quad
D^{(1)}_\zs{2,n}\,
\le 
 \phi^{2}_\zs{max}\,
\vert\rho\vert_\zs{*}\,n^{2}
\,.
\end{equation}
Let us estimate now the last term in \eqref{sec:Pr.35-00-1}. First, note that we can represent this term as
$$
D^{(2)}_\zs{2,n}=\sum_\zs{k=1}^{+\infty} \E\,
\Chi_\zs{\{T_\zs{k} \leq n\}}\,
\sum^{k-1}_\zs{j=1}\,
Q^{2}_\zs{x}(T_\zs{k},T_\zs{j})=
\sum^{\infty}_\zs{j=1}
\,\Chi_\zs{\{T_\zs{j} \leq n\}}\,G(T_\zs{j})
=\int^{n}_\zs{0}\,G(t)\,\rho(t)\d t\,,
$$
where
\begin{align*}
G(t)=&
\sum_\zs{k=1}^{+\infty} \E\,
\Chi_\zs{\{T_\zs{k} \leq n\}}\,
Q^{2}_\zs{x}((t+T_\zs{k}),t)
=
\int^{n}_\zs{0}\,
Q^{2}_\zs{x}(t+v,t)
\,
\rho(v)\d v\\[2mm]
=&
\int^{n+t}_\zs{t}\,
Q^{2}_\zs{x}(u,t)
\,
\rho(u-t)\d u
\,.
\end{align*}
It is clear that, for any $0\le t\le n,$
$$
\int^{n+t}_\zs{t}\,
Q^{2}_\zs{x}(u,u-t)
\,
\rho(u)\,
\d u\,
\le \vert\rho\vert_\zs{*}\,
\int^{2n}_\zs{0}\,
Q^{2}_\zs{x}(v,t)
\,\d v\,.
$$
In view of the inequality \eqref{sec:Pr.35-00-1-2} we obtain
$$
\int^{2n}_\zs{0}\,
Q^{2}_\zs{x}(u,t)
\,\d u
=
\int^{2n}_\zs{0}\,
Q^{2}_\zs{x}(t,u)
\,\d u
\,
\le 
2 \phi^{2}_\zs{max}\,n\,.
$$
Therefore,
$$
\max_\zs{0\le t\le n}\,G(t)\,
\le\,
2 \vert\rho\vert_\zs{*}\,\phi^{2}_\zs{max}\,n
\quad\mbox{and}\quad
D^{(2)}_\zs{2,n}\le
\,
2 \vert\rho\vert^{2}_\zs{*}\,\phi^{2}_\zs{max}\,
n^{2}\,.
$$
So, estimating  $\varrho^{2}_\zs{2}$  by $\check{\tau}\sigma_\zs{Q}$ 
and taking into account that $\E Y^{4}_\zs{1}\ge 1,$ we obtain that
we obtain that
$$
\E \,D_n(x) \le   
13\,(1+\check{\tau})
\phi^{4}_\zs{max}\,
\,
\E Y_1^4 
(1+\vert\rho\vert^{2}_\zs{*})\,
n^{2}
\sigma_\zs{Q}
\,.
$$
Using all these bound in \eqref{sec:Pr.4-11-03} we obtain \eqref{sec:OI.1-14.3.1-2} and thus the conclusion follows. \fdem

\begin{remark}\label{Re.sec:Prsm.1}
The properties  \eqref{sec:OI.1-14.3.1} and \eqref{sec:OI.1-14.3.1-2} are used to obtain
the oracle inequalities given in Section \ref{sec:Mrs}
 (see, for example, \cite{KonevPergamenshchikov2012}).
\end{remark}

\section{Simulation}\label{sec:Siml}

In this section we report the results of a Monte Carlo experiment in order to assess the performance of the proposed model selection procedure \eqref{sec:Mo.9}. 
In \eqref{sec:In.1} we chose a  $1$-periodic function which is defined, for $0\le t\le 1,$ as 
\begin{equation}\label{sec:Siml.0}
S(t)=t \sin(2\pi t)+t^2(1-t) \cos(4\pi t)\,.
\end{equation}
We simulate the model
$$
\d y_\zs{t} = S(t) \d t + \d\xi_\zs{t}\,, 
$$
where $\xi_t= 0.5 \d w_\zs{}t+ 0.5 \d z_\zs{t}$. 

Here $z_\zs{t}$ is the semi-Markov process defined in \eqref{sec:Ex.2} with a Gaussian $\mathcal{N}(0,1)$ 
sequence $(Y_\zs{j})_\zs{j\geq1}$ and $(\tau_k)_{k\geq1}$ used in \eqref{sec:Ex.4} taken as
$\tau_k \sim \chi_\zs{3}^2.$ 

We use the model selection procedure \eqref{sec:Mo.9} with the weights \eqref{sec:Ga.2}
in which $k^*= 100+\sqrt{\ln(n)}$, $t_\zs{i}=i/ \ln (n)$, $m=[\ln^2 (n)]$ and $\delta=(3+\ln(n))^{-2}$.
We define the empirical risk as
\begin{equation}\label{sec:Siml.1}
\mathbf{\overline{R}}= \frac{1}{p} \sum_{j=1}^{p}  \mathbf{\hat{E}} \left(\wh{S}_n(t_j)-S(t_j)\right)^2,
\end{equation}
where the observation frequency   $p=100001$ and  the expectation was taken as an average over $N= 10000$ replications, i.e.,
$$
\mathbf{\hat{E}} \left(\wh{S}_n(.)-S(.)\right)^2 = \frac{1}{N} \sum_{l=1}^{N} \left(\wh{S}^l_n(.)-S(.) \right)^2.
$$
\noindent 
We set the relative quadratic risk as
\begin{equation}\label{sec:Siml.2}
\mathbf{\overline{R_*}}=\mathbf{\overline{R}}/ ||S||^2_\zs{p},
\quad\mbox{with}\quad
||S||^2_p = \frac{1}{p} \sum_{j=0}^p S^2(t_j)\,.
\end{equation}
In our case $||S||^2_p = 0.1883601$.

Table \ref{tab:1} gives the values for the sample risks \eqref{sec:Siml.1} and \eqref{sec:Siml.2} for different numbers of observations $n$.

{\renewcommand{\arraystretch}{2} 
{\setlength{\tabcolsep}{1cm} 

\begin{table}  
\begin{center}
    \begin{tabular}{|r|c|c|}
                                  \hline
                                  n & $\mathbf{\overline{R}}$  & $\mathbf{\overline{R_*}}$  \\
                                  \hline
                                  20 &0.04430  & 0.235  \\
                                  \hline
                                  100 & 0.01290 &0.068  \\
                                  \hline
                                  200 & 0.00812 &0.043  \\
                                  \hline
                                  1000 &0.00196  &0.010  \\
                                  \hline
		\end{tabular}
\end{center}
\caption{Empirical risks} \label{tab:1}
\end{table}

Figures \ref{fig1}--\ref{fig4} show  the behaviour of the regression function and its estimates
by the model selection procedure \eqref{sec:Mo.9}
depending on the values of observation periods $n$. The black full line is the regression function \eqref{sec:Siml.0}
and the red dotted line is the associated  estimator.


\newpage

\begin{figure}[!h]
\begin{center}
\includegraphics[scale=0.25]{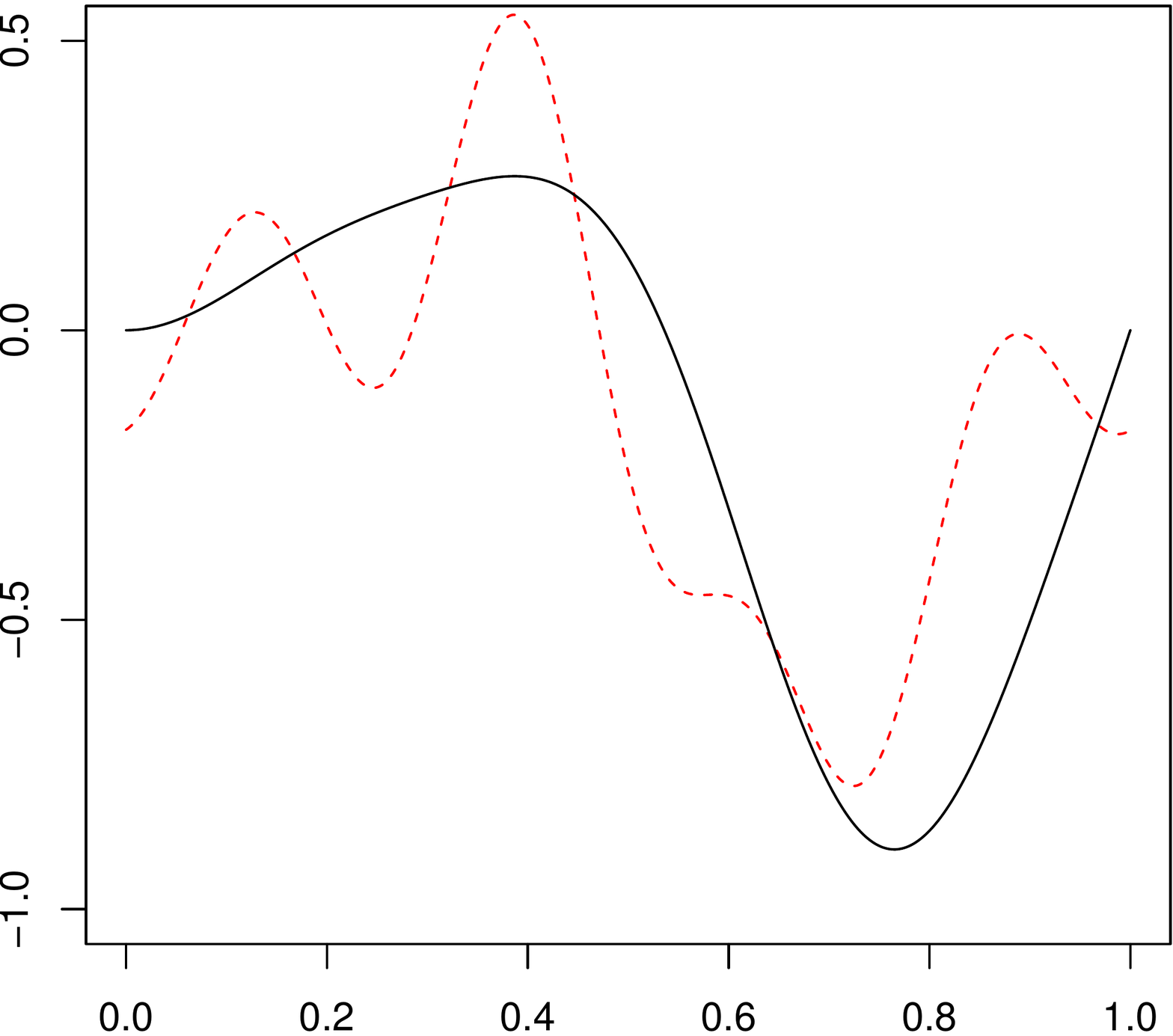}
\caption{Estimator of $S$ for $n=20$} \label{fig1}
\includegraphics[scale=0.25]{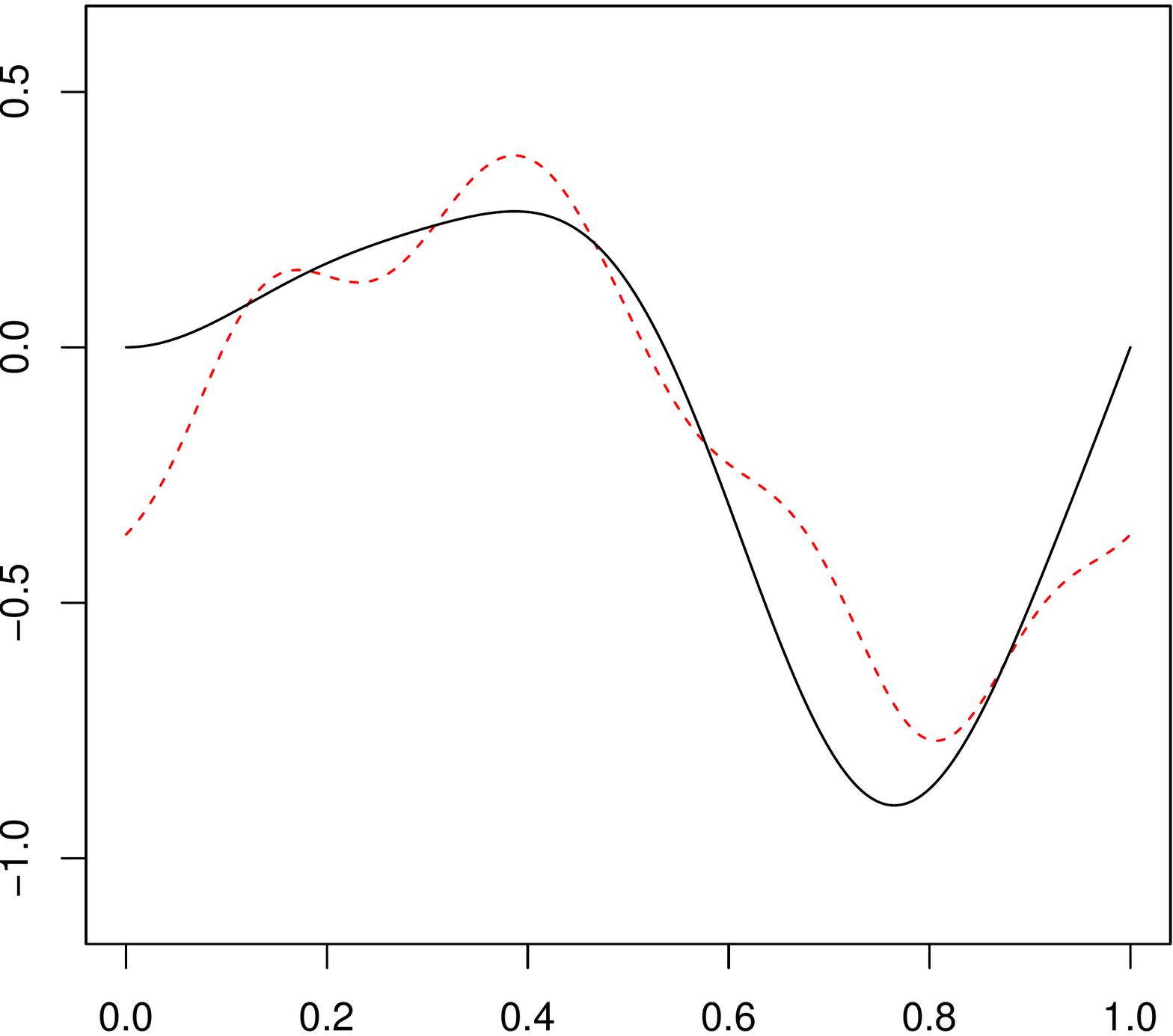}
\caption{Estimator of $S$ for $n=100$} \label{fig2}

\end{center}
\end{figure}

\newpage



\begin{figure}[!h]
\begin{center}
\includegraphics[scale=0.25]{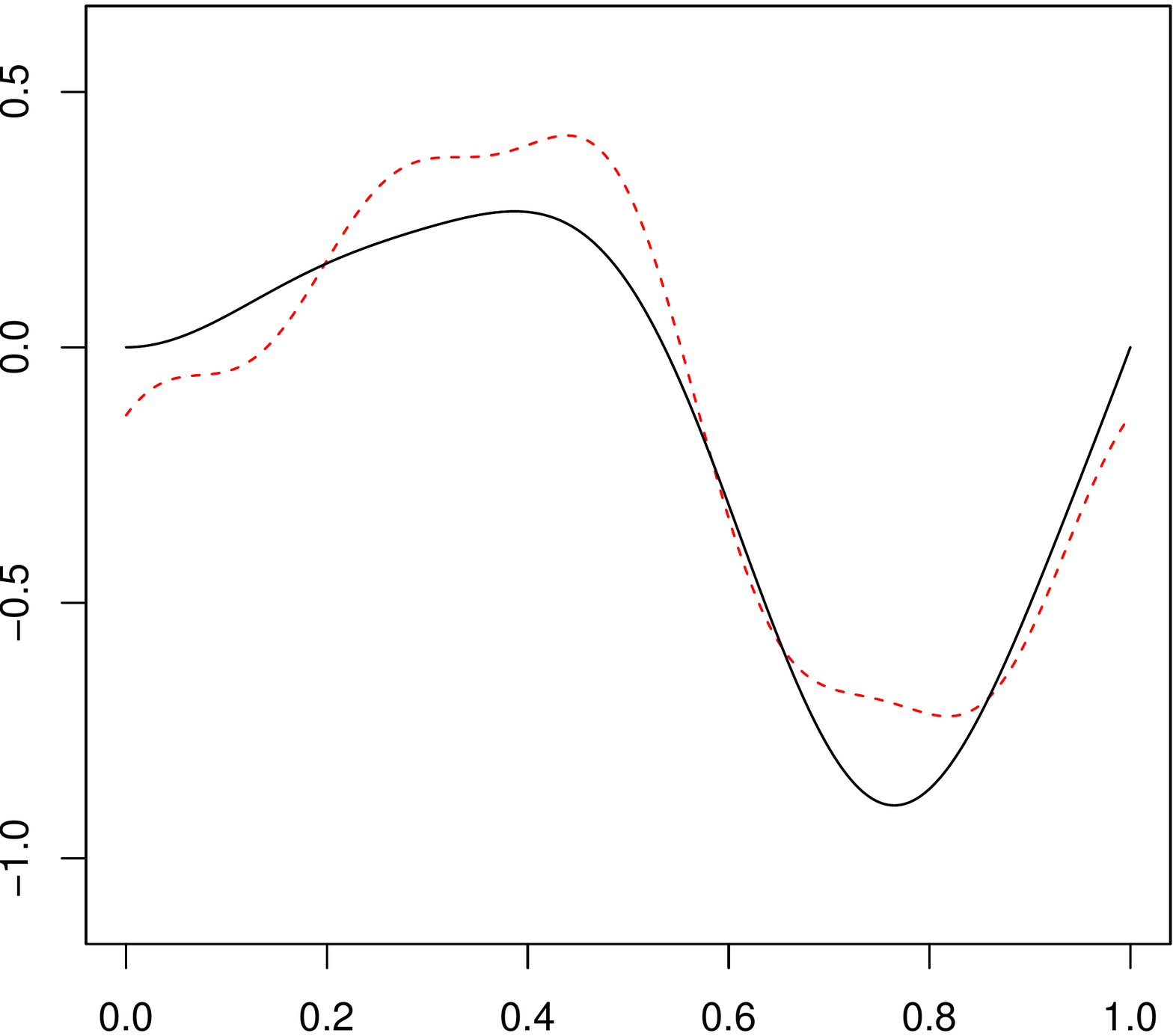}
\caption{Estimator of $S$ for $n=200$}\label{fig3}
\includegraphics[scale=0.25]{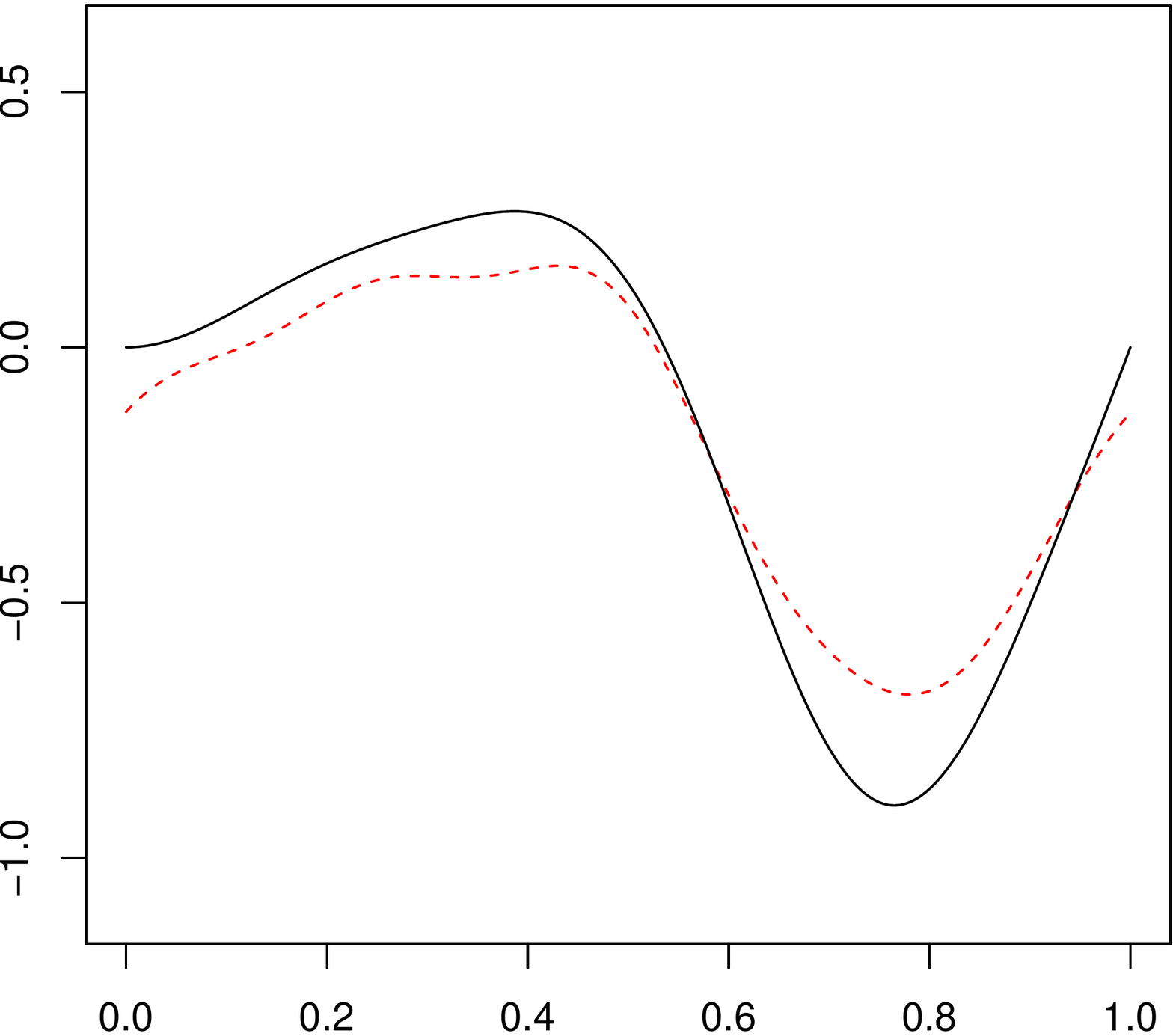}
\caption{Estimator of $S$ for $n=1000$} \label{fig4}
\end{center}
\end{figure}

\newpage

\begin{remark}\label{Re.sec:Mnc.1}
From numerical simulations  of the procedure \eqref{sec:Mo.9}
with various observation numbers $n$ we may conclude that
 the quality
of the proposed procedure: (i) is good for practical needs, i.e. for reasonable (non large)  number of observations; (ii) is improving as the number of observations increases.
\end{remark}

\section{Proofs}\label{sec:Pr}

We will prove here most of the results of this paper.

\subsection{Proof of Theorem \ref{Th.sec:OI.1}} 

First, note that
 we can rewrite the empirical squared error in \eqref{sec:Mo.3} as follows
\begin{equation}\label{sec:Pr.1}
\Er_n(\lambda) = J_n(\lambda) + 2 \sum_\zs{j=1}^{\infty} \lambda(j) \check{\theta}_\zs{j,n}+ ||S||^2-\delta P_n(\lambda),
\end{equation}
where $\check{\theta}_\zs{j,n}=\wt{\theta}_\zs{j,n}-\theta_\zs{j}\wh{\theta}_\zs{j,n}$. Using the definition of $\wt{\theta}_\zs{j,n}$ in \eqref{sec:Mo.4}
we obtain that
$$
\check{\theta}_\zs{j,n}=\frac{1}{\sqrt n}\theta_\zs{j}\xi_\zs{j,n} +\frac{1}{n}\wt{\xi}_\zs{j,n}  + \frac{1}{n} \varsigma_\zs{j,n}  +
\frac{\sigma_\zs{Q} -  \wh{\sigma}_\zs{n} }{n}\,,
$$
where
$\varsigma_\zs{j,n}=\E_\zs{Q}\xi^{2}_\zs{j,n}-\sigma_\zs{Q}$ and $\wt{\xi}_\zs{j,n}=\xi^{2}_\zs{j,n}-\E_\zs{Q}\xi^{2}_\zs{j,n}$. Putting
\begin{equation}\label{sec:Pr.2}
M(\lambda) = \frac{1}{\sqrt n}\sum_\zs{j=1}^{n} \lambda(j)\theta_\zs{j} \xi_\zs{j,n}
\quad\mbox{and}\quad
P^{0}_\zs{n}=\frac{\sigma_\zs{Q}\vert\lambda\vert^{2}}{n}\,,
\end{equation}
we can rewrite \eqref{sec:Pr.1} as
\begin{align}\nonumber
\Er_n(\lambda)  = &  J_n(\lambda) + 2 \frac{\sigma_\zs{Q}-  \wh{\sigma}_\zs{n} }{n}\,\check{L}(\lambda)+ 2 M(\lambda)+\frac{2}{n} B_\zs{1,Q,n}(\lambda)\\  \label{sec:Pr.3}
& +  2 \sqrt{P^{0}_n(\lambda)} \frac{B_\zs{2,Q,n}(e(\lambda)}{\sqrt{\sigma_\zs{Q} n}} + \Vert S\Vert^2-\rho  P_n(\lambda),
\end{align}
where $e(\lambda)=\lambda/|\lambda|$, the function $\check{L}(\cdot)$ is defined in \eqref{sec:Mo.2} and the functions $ B_\zs{1,Q,n}(\cdot)$ and $B_\zs{2,Q,n}(\cdot)$
 are given in \eqref{sec:Prsm.1}.

Let $\lambda_0= (\lambda_0(j))_\zs{1\le j\le\,n}$ be a fixed sequence in $\Lambda$ and $\wh{\lambda}$ be as in \eqref{sec:Mo.8}.
Substituting $\lambda_0$ and $\wh{\lambda}$ in Equation \eqref{sec:Pr.3}, we obtain
\begin{align}\label{sec:Pr.4}
\Er_n(\wh{\lambda})-\Er_n(\lambda_0) = & J(\wh{\lambda})-J(\lambda_0)+
2 \frac{\sigma_\zs{Q}-\wh{\sigma}_\zs{Q}}{n}\,\check{L}(\varpi)
+ \frac{2}{n} B_\zs{1,Q,n}(\varpi)+2 M(\varpi)\nonumber\\[2mm]
& + 2 \sqrt{P^{0}_\zs{n}(\wh{\lambda})} \frac{B_\zs{2,Q,n}(\wh e)}{\sqrt{\sigma_\zs{Q} n}}-2 \sqrt{P^{0}_\zs{n}(\lambda_0)}
 \frac{B_\zs{2,Q,n}(e_0)}{\sqrt{\sigma_\zs{Q} n}}\nonumber \\[2mm]
& -  \delta  P_n(\wh{\lambda})+\delta P_n(\lambda_0),
\end{align}
where $\varpi= \wh{\lambda} - \lambda_\zs{0}$, $\wh{e} = e(\wh{\lambda})$ and $e_0 = e(\lambda_0)$. Note that, by \eqref{sec:Mo.2},
$$ 
|\check{L}(\hat x)| \le\,\check{L}(\hat \lambda) + \check{L}(\lambda) \leq 2\vert\Lambda\vert_\zs{*}. 
$$
Applying the inequality
\begin{equation}\label{sec:Pr.5}
2|ab| \leq \delta a^2 + \delta^{-1} b^2
\end{equation}
implies that, for any $\lambda\in\Lambda,$
$$
2 \sqrt{P^{0}_n(\lambda)} \frac{|B_\zs{2,Q,n}(e(\lambda))|}{\sqrt{\sigma_\zs{Q} n}} \le\, \delta P^{0}_\zs{n}(\lambda) +
 \frac{B^2_\zs{2,Q,n}(e(\lambda))}{\delta\sigma_\zs{Q}\,n}.
$$
Taking into account the bound \eqref{sec:OI.1-14.3.1}, we get
\begin{align*}
\Er_n(\hat \lambda)  \le &\,\Er_n(\lambda_0) +2 M(\varpi)+ \frac{2 \C_\zs{1,Q,n}}{n}+ \frac{2 B^*_\zs{2,Q,n}}{\delta\sigma_\zs{Q}\,n} \\[2mm]
& +  \frac{1}{n} |\wh{\sigma} -\sigma_\zs{Q}| ( |\wh{\lambda}|^2 + |\lambda_0|^2)+ 2  \delta P_n(\lambda_0)\,,
\end{align*}
where $B^*_\zs{2,Q,n} = \sup_\zs{\lambda\in\Lambda} B^2_\zs{2,Q,n}((e(\lambda))$. 
Moreover, noting that in view of \eqref{sec:Mo.2} $\sup_\zs{\lambda\in\Lambda} |\lambda|^2 \leq \vert\Lambda\vert_\zs{*}$,
we can rewrite the previous  bound as
\begin{align}\label{sec:Pr.6}
\Er_n(\wh{\lambda})  \le & \Er_n(\lambda_0) +2 M(\varpi)
+ \frac{2 \C_\zs{1,Q,n}}{n}
+ \frac{2 B^*_\zs{2,Q,n}}{\delta\sigma_\zs{Q} n}  \nonumber\\[2mm]
& +  \frac{4\vert\Lambda\vert_\zs{*}}{n} |\wh{\sigma} -\sigma_\zs{Q}| + 2  \delta P_n(\lambda_0).
\end{align}
To estimate the second term in the right side of this inequality we set
$$
S_\zs{x} = \sum_\zs{j=1}^{n} x(j) \theta_\zs{j} \phi_\zs{j}\,,
\quad x=(x(j))_\zs{1\le j\le n}\in\bbr^{n}\,.
$$
Thanks to 
\eqref{sec:In.3} we estimate the term $M(x)$ for any $x\in\bbr^{n}$ as
\begin{equation}\label{sec:Pr.8}
\E_\zs{Q} M^2 (x) \leq \varkappa_\zs{Q} \frac{1}{n} \sum_\zs{j=1}^{n} x^2(j) \theta^2_\zs{j} = \varkappa_\zs{Q}\frac{1}{n} ||S_\zs{x}||^2.
\end{equation}
To estimate this function for a random vector $x\in\bbr^{n}$ we set
$$ 
Z^* = \sup_\zs{x \varepsilon \Lambda_1} \frac{n M^2 (x)}{||S_x||^2}\,,
\quad
\Lambda_1 = \Lambda - \lambda_0\,.
$$
So, through Inequality  \eqref{sec:Pr.5}, we get
\begin{equation}\label{sec:Pr.10}
2 |M(x)|\leq \delta ||S_x||^2 + \frac{Z^*}{n\delta}.
\end{equation}
It is clear that the last term  here can be estimated as
\begin{equation}\label{sec:Pr.9}
\E_\zs{Q} Z^* \leq \sum_\zs{x \in \Lambda_1} \frac{n \E_\zs{Q} M^2 (x)}{||S_x||^2} \leq \sum_\zs{x \in \Lambda_1} \varkappa_\zs{Q}= \varkappa_\zs{Q}\check{\iota}\,,
\end{equation}
where $\check{\iota} = \mbox{card}(\Lambda)$. 
Moreover, note that, for any $x\in\Lambda_1$,
\begin{equation}\label{sec:Pr.11}
||S_x||^2-||\wh{S}_x||^2 = \sum_\zs{j=1}^{n} x^2(j) (\theta^2_\zs{j}-\wh{\theta}^2_\zs{j}) \le -2 M_1(x),
\end{equation}
where $M_\zs{1}(x) =  n^{-1/2}\,\sum_\zs{j=1}^{n}\, x^2(j)\theta_\zs{j} \xi_\zs{j,n}$.
Taking into account that, for any $x \in \Lambda_1$ the components $|x(j)|\leq 1$ ,  we can estimate this term as
in \eqref{sec:Pr.8}, i.e.,
$$
\E_\zs{Q}\, M^2_\zs{1}(x) \leq \varkappa_\zs{Q}\,
\frac{||S_x||^2}{n}\,.
$$
Similarly to the previous reasoning
we set
$$ 
Z^*_\zs{1} = \sup_\zs{x \varepsilon \Lambda_1} \frac{n M^2_1 (x)}{||S_x||^2}
$$
and we get
\begin{equation}\label{sec:Pr.12}
\E_\zs{Q}\, Z^*_1 \leq \varkappa_\zs{Q}\,\check{\iota}\,.
\end{equation}
Using the same type of arguments as in \eqref{sec:Pr.10}, we can derive
\begin{equation}\label{sec:Pr.13}
2 |M_1(x)|\leq \delta ||S_x||^2 + \frac{Z^*_1}{n\delta}.
\end{equation}
From here and \eqref{sec:Pr.11}, we get
\begin{equation}\label{sec:Pr.14}
||S_x||^2 \leq \frac{||\wh{S}_x||^2}{1-\delta} + \frac{Z^*_1}{n \delta (1-\delta)}
\end{equation}
for any $0<\delta<1$. Using this bound in \eqref{sec:Pr.10} yields
$$
2 M(x) \leq \frac{\delta ||\wh{S}_x||^2}{1-\delta} + \frac{Z^*+Z^*_1}{n \delta (1-\delta)}
\,.
$$
Taking into account that $\Vert\wh{S}_\zs{\varpi}\Vert^{2}\le 2\,(\Er_n(\wh{\lambda})+\Er_n(\lambda_0))$, we obtain
$$
2 M(\varpi) \leq \frac{2\delta(\Er_n(\wh{\lambda})+\Er_n(\lambda_0))}{1-\delta} + \frac{Z^*+Z^*_1}{n \delta (1-\delta)}.
$$
Using this bound in  \eqref{sec:Pr.6} we obtain 
\begin{align*}
\Er_n(\wh{\lambda})  \le & \frac{1+\delta}{1-3\delta} \Er_n(\lambda_0) 
+ \frac{Z^*+Z^*_1}{n \delta (1-3\delta)}
+ \frac{2 \C_\zs{1,Q,n}}{n(1-3\delta)}
+ \frac{2 B^*_\zs{2,Q,n}}{\delta(1-3\delta)\sigma_\zs{Q} n} \\[2mm]
& +  \frac{(4\vert\Lambda\vert_\zs{*}+2)}{n(1-3\delta)} |\wh{\sigma} -\sigma_\zs{Q}| + \frac{2\delta}{(1-3\delta)} P^{0}_n(\lambda_0).
\end{align*}
Moreover, for $0<\delta<1/6,$ we can rewrite this inequality as
\begin{align*}
\Er_n(\wh{\lambda})  \le & \frac{1+\delta}{1-3\delta} \Er_n(\lambda_0) 
+ \frac{2(Z^*+Z^*_1)}{n \delta}
+ \frac{4 \C_\zs{1,Q,n}}{n}
+ \frac{4 B^*_\zs{2,Q,n}}{\delta \sigma_\zs{Q} n} \\[2mm]
& +  \frac{(8\vert\Lambda\vert_\zs{*}+2)}{n} |\wh{\sigma}_\zs{n} -\sigma_\zs{Q}| + \frac{2\delta}{(1-3\delta)}\,
 P^{0}_n(\lambda_0).
\end{align*}
In view of Proposition \ref{Pr.sec:OI.3} we estimate the expectation of the term $B^*_\zs{2,Q,n}$ in \eqref{sec:Pr.6} as
$$
\E_\zs{Q}\, B^*_\zs{2,Q,n} \leq \sum_\zs{\lambda\in\Lambda}\E_\zs{Q} B^2_\zs{2,Q,n} (e(\lambda)) \leq \check{\iota} \C_\zs{2,Q,n}\,.
$$
Taking into account that $\vert\Lambda\vert_\zs{*}\ge 1$, we get
\begin{align*}
 \cR(\wh{S}_*,S)  \le & \frac{1+\delta}{1-3\delta} 
 \cR(\wh{S}_\zs{\lambda_0},S)
+ \frac{4\varkappa_\zs{Q} \check{\iota}}{n \delta}
+ \frac{4 \C_\zs{1,Q,n}}{n}
+ \frac{4 \check{\iota} \C_\zs{2,Q,n}}{\delta \sigma_\zs{Q} n} \\[2mm]
& +  \frac{10\vert\Lambda\vert_\zs{*}}{n} \,\E_\zs{Q}\,|\wh{\sigma} -\sigma_\zs{Q}| + \frac{2\delta}{(1-3\delta)} P^{0}_n(\lambda_0).
\end{align*}
Using the upper bound for $ P_n(\lambda_0)$ in Lemma~\ref{Le.sec:A.1-06-11-01}, one obtains 
\eqref{Th.sec:OI.1}, that finishes the proof. \fdem

\subsection{Proof of Proposition \ref{Pr.sec:Si.1}}
We use here the same method as in \cite{KonevPergamenshchikov2009a}.
First of all note that the definition \eqref{sec:Mo.4-2-31-3} implies that
\begin{equation}\label{sec:Mo.1-1-04}
\wh{t}_\zs{j,n}= t_\zs{j}+
\frac{1}{\sqrt{n}}\,
\eta_\zs{j,n}\,,
\end{equation}
where  
$$
t_\zs{j}=
\int^{1}_\zs{0}\,S(t)\,Tr_\zs{j}(t)\d t
\quad\mbox{and}\quad
\eta_\zs{j,n}=
\frac{1}{\sqrt{n}}\,
\int^{n}_\zs{0}\,\Tr_\zs{j}(t)\,\d \xi_\zs{t}
\,.
$$
So, we have
\begin{equation}\label{sec:Mo.1-0-1-04}
\wh{\sigma}_\zs{n}=
\sum^n_\zs{j=[\sqrt{n}]+1}\,t^2_\zs{j}
+
2M_\zs{n}
+
\frac{1}{n}\,
\sum^n_\zs{j=[\sqrt{n}]+1}\,\eta^2_\zs{j,n}
\,,
\end{equation}
where 
$$
M_\zs{n}=
\frac{1}{\sqrt{n}}
\sum^n_\zs{j=[\sqrt{n}]+1}\,t_\zs{j}\,\eta_\zs{j,n}\,.
$$
Note that, for continuously differentiable functions (see, for example, Lemma A.6 in \cite{KonevPergamenshchikov2009a}),
the Fourier coefficients $(t_\zs{j})$ satisfy the following inequality, for any $n\ge 1,$
\begin{equation}\label{sec:Mo.2-1-04}
\sum^{\infty}_\zs{j=[\sqrt{n}]+1}\,t^2_\zs{j}
\le 
\frac{4\left(\int^{1}_\zs{0}\vert\dot{S}(t)\vert \d t\right)^{2}}{\sqrt{n}}
\le 
\frac{4\Vert\dot{S}\Vert^{2}}{\sqrt{n}}
\,.
\end{equation}
In the same way as in
\eqref{sec:Pr.8} we estimate the term $M_\zs{n}$, i.e.,
$$
\E_\zs{Q}\,M^{2}_\zs{n}\le \frac{\varkappa_\zs{Q}}{n}\,
\sum^{n}_\zs{j=[\sqrt{n}]+1}\,t^{2}_\zs{j}
\le 
\frac{4\varkappa_\zs{Q}\Vert\dot{S}\Vert^{2}}{n\sqrt{n}},
$$
while the absolute value of this term for $n\ge 1$ can be estimated  as
$$
\vert \E_\zs{Q}\,M_\zs{n}\vert
\le 
\frac{\varkappa_\zs{Q}+\Vert\dot{S}\Vert^{2}}{\sqrt{n}}\,.
$$
Moreover, using Propositions  \ref{Pr.sec:OI.2} and \ref{Pr.sec:OI.3} we can represent the last term in \eqref{sec:Mo.1-0-1-04} as
$$
\frac{1}{n}
\sum^n_\zs{j=[\sqrt{n}]+1}\,\eta^2_\zs{j,n}
=\frac{\sigma_\zs{Q}(n-\sqrt{n})}{n}
+\frac{B_\zs{1,Q,n}(x')}{n}
+\frac{B_\zs{2,Q,n}(x'')}{\sqrt{n}},
$$
with
$
x'_\zs{j}=\Chi_\zs{\{\sqrt{n}<j\le n\}}$
and $x''_\zs{j}=\Chi_\zs{\{\sqrt{n}<j\le n\}}/\sqrt{n}$. 
Therefore, 
$$
\E_\zs{Q}
\left\vert
\frac{1}{n}
\sum^n_\zs{j=[\sqrt{n}]+1}\,\eta^2_\zs{j,n}
-\sigma_\zs{Q}
\right\vert
\le 
\frac{\sigma_\zs{Q}}{\sqrt{n}}
+\frac{\C_\zs{1,Q,n}}{n}
+\frac{\sqrt{\C_\zs{2,Q,n}}}{\sqrt{n}}.
$$
Taking into account that $\C_\zs{2,Q,n}\ge 1,$
we obtain the bound \eqref{sec:Si.3} and hence the desired result.
\fdem 
\subsection{Proof of Theorem \ref{Th.sec:Mrs.2}}

First note, that  in view of  \eqref{sec:Ga.1++1--1}
and  \eqref{sec:Ga.1}
$$
\lim_\zs{n\to\infty}\,\frac{\check{\iota}}{n^{\check{\epsilon}}}=
\lim_\zs{n\to\infty}\,\frac{k^{*} m}{n^{\check{\epsilon}}}=0
\qquad\mbox{for any}\quad 
\check{\epsilon}>0
\,.
$$
Furthermore, the bound \eqref{sec:Ga.1++1--2} and the  conditions \eqref{sec:Mrs.5-1} and \eqref{sec:Ga.1} yield
$$
\lim_\zs{n\to\infty}\frac{\vert\Lambda\vert_\zs{*}}{n^{1/3+\check{\epsilon}}}=0
\quad\mbox{for any}\quad \check{\epsilon}>0
\,.
$$
So, from here we obtain the convergence \eqref{sec:Mrs.7-25.3}. \fdem

 \bigskip
 
 \subsection{Proof of Theorem \ref{Th.sec:Ef.1}}
 
 First, we denote by $Q_\zs{0}$ the distribution of the noise
 \eqref{sec:Ex.1} and \eqref{sec:Mcs.1}
 with the parameter $\varrho_\zs{1}=\varsigma^{*}$, $\check{\varrho}=1$ and $\varrho_\zs{2}=0$, i.e.
 the distribution
 for the ``signal + white noise''  
 model. So, we can estimate as below the robust risk
 $$
\cR^{*}_\zs{n}(\wt{S}_\zs{n},S)\ge
\cR_\zs{Q_\zs{0}}(\wt{S}_\zs{n},S)\,.
$$
Now Theorem 6.1 from \cite{KonevPergamenshchikov2009b}
yields the lower bound \eqref{sec:Ef.4}.  Hence this finishes the proof. \fdem

 \subsection{Proof of Proposition \ref{Th.sec:Ef.33}}
 
Putting $\lambda_\zs{0}(j)=0$ for $j\ge n$ we can represent the quadratic risk for the estimator \eqref{sec:Mo.1} as 
$$
\parallel \wh{S}_\zs{\lambda_0}-S\parallel^2=\sum_\zs{j=1}^{\infty} (1-\lambda_0(j))^2 \theta^2_\zs{j}-2 H_\zs{n}
+ \frac{1}{n} \sum_{j=1}^{n}   \lambda_0^2(j) \xi^2_{j,n}\,, 
$$
where $H_n= n^{-1/2}\,\sum_{j=1}^{n} (1-\lambda_0(j)) \lambda_0(j) \theta_{j} \xi_\zs{j,n}$. 
Note that $\E_{Q} H_\zs{n}=0$ for any  $ Q \in Q_n$, therefore,
$$
\E_\zs{Q}\parallel \wh{S}_\zs{\lambda_0}-S\parallel^2=\sum_\zs{j=1}^{\infty} (1-\lambda_0(j))^2 \theta^2_\zs{j}
+ \frac{1}{n} \E_\zs{Q}\sum_{j=1}^{n}   \lambda_0^2(j) \xi^2_{j,n}\,.
$$
Proposition \ref{Pr.sec:OI.2} and the last inequality in \eqref{sec:Ex.5}
imply that for any $Q\in\cQ_\zs{n}$
$$ 
\E_{Q} \sum_{j=1}^{n} \lambda_0^2(j) \xi^2_{j,n} \leq  
\varsigma^*   \sum_{j=1}^{n} \lambda_0^2(j) + 
\frac{\phi^{2}_\zs{max}\varsigma^{*}\Vert\Upsilon\Vert_\zs{1}}{\check{\tau}}
:=
\varsigma^*   \sum_{j=1}^{n} \lambda_0^2(j) +
\C^{*}_\zs{1,n}\,.           
$$
Therefore,
$$
\cR^*_\zs{n} (\wh{S}_{\lambda_\zs{0}},S) \leq \sum_\zs{j=j_\zs{*}}^{\infty} (1-\lambda_0(j))^2 \theta^2_{j}+
\frac{1}{\upsilon_\zs{n}}  \sum_{j=1}^{n}   \lambda_0^2(j)  
+\frac{\C^{*}_\zs{1,n}}{n},
$$
where $j_\zs{*}$ and $\upsilon_\zs{n}$ are defined in \eqref{sec:Ga.2}. Setting
$$
\Upsilon_\zs{1,n}(S) = \upsilon^{2k/(2k+1)}_\zs{n} \sum_\zs{j=j_\zs{*}}^{\infty} (1-\lambda_0(j))^2 \theta^2_{j}
\quad\mbox{and}\quad
\Upsilon_\zs{2,n}=
\frac{1}{\upsilon^{1/(2k+1)}_\zs{n}}  \sum_{j=1}^{n}   \lambda_0^2(j) \,,
$$
we rewrite the last inequality as
\begin{equation}\label{sec:Ef.2.1++}
\upsilon^{2k/(2k+1)}_\zs{n}\,
R^*_\zs{n} (\wh{S}_{\lambda_\zs{0}},S) \leq 
\Upsilon_\zs{1,n}(S)
+
\Upsilon_\zs{2,n}
+
\check{\C}_\zs{n}
\,,
\end{equation}
where $
\check{\C}_\zs{n}=
\upsilon^{2k/(2k+1)}_\zs{n}\C^{*}_\zs{1,n}/n$. 
Note, that the conditions \eqref{sec:Mrs.5-1} and \eqref{sec:Mrs.5-2}
imply  that $\C^{*}_\zs{1,n}=\oo(n^{\check{\delta}})$  as $n\to\infty$ for any $\check{\delta}>0$; therefore, 
$\check{\C}_\zs{n}\to 0$ as $n\to\infty$.
Putting
$$
u_\zs{n}= \upsilon^{2k/(2k+1)}_\zs{n} \sup_\zs{j\geq j_\zs{*}} (1-\lambda_0(j))^2/a_j
\,,
$$
with $a_j$ defined in \eqref{sec:Ef.2},
we estimate the first term in \eqref{sec:Ef.2.1++} as 
$$
\sup_\zs{S\in W_\zs{\r}^k}\,
\Upsilon_\zs{1,n}(S)
\le
\sup_\zs{S\in W_\zs{\r}^k}\, u_\zs{n}\,\sum_{j\geq1} a_j \theta_{j}
 \le 
 u_\zs{n} \r
 \,.
$$
Taking into account that $a_\zs{j}/(\pi^{2k}j^{2k})\to 1$ as $j\to \infty$ and $\l_\zs{0}\to \r$ as $\varepsilon\to 0$
and using the definition of $\omega_\zs{\alpha_\zs{0}}$ in \eqref{sec:Ga.2},
we obtain that
\begin{align*}
\limsup_\zs{n\to \infty }  u_\zs{n} &\leq 
\lim_\zs{n\to \infty } 
\,
\upsilon^{2k/(2k+1)}_\zs{n} \sup_\zs{j\geq j_\zs{*}} \frac{(1-\lambda_0(j))^2}{(\pi\,j)^{2k}}\\[2mm]
&=
\lim_\zs{n\to \infty } 
\,
\frac{\upsilon^{2k/(2k+1)}_\zs{n}}{\pi^{2k} \omega^{2k}_\zs{\alpha_\zs{0}}}
=
\frac{1}{\pi^{2k}\,(\d_\zs{k}\r)^{2k/(2k+1)}}\,.
\end{align*}
Therefore,
\begin{equation}
\limsup_\zs{n\to\infty } \sup_\zs{S\in W^{k}_\zs{\r}} \Upsilon_{1,n}(S) \leq
 \frac{r^{1/(2k+1)}}{\pi^{2k} (\d_\zs{k})^{2k/(2k+1)}}=:
 \Upsilon^*_\zs{1}\,.
\end{equation}
As to the second term in \eqref{sec:Ef.2.1++}, note that
$$
\lim_\zs{n\to\infty} \frac{1}{\omega_\zs{\alpha_\zs{0}}}\,\sum_{j=1}^{n}   \lambda_0^2(j)=
\int^{1}_\zs{0}(1-t^{k})^{2}\d t=
\frac{2k^{2}}{(k+1)(2k+1)}
\,.
$$
So, taking into account that $\omega_\zs{\alpha_\zs{0}}/\upsilon^{1/(2k+1)}_\zs{n}\to ( \d_\zs{k} \r)^{1/(2k+1)}$ as $n\to\infty$, the limit
of $\Upsilon_\zs{2,n}$ can be calculated as
$$
\lim_\zs{n\to \infty } \,\Upsilon_\zs{2,n}= 
\frac{2( \d_\zs{k} \r)^{1/(2k+1)}\,k^{2}}{(k+1)(2k+1)}=: \Upsilon^*_\zs{2}\,.    
$$
Moreover, since $\Upsilon^*_1+ \Upsilon^*_2 =:\r^*_\zs{k}$,
we obtain 
$$
\lim_\zs{n \to \infty } \upsilon^{2k /(2k+1)}_\zs{n}\, \sup_\zs{S\in W^{k}_\zs{\r}} \cR^*_n (\wh{S}_{\lambda_\zs{0}},S) \leq \r^*_\zs{k}
$$
and get the desired result.  \fdem

\subsection{Proof of Theorem \ref{Th.sec:Ef.2}}
Combining Proposition \ref{Th.sec:Ef.33} and  Theorem \ref{Th.sec:Mrs.2} yields Theorem  \ref{Th.sec:Ef.33}. \fdem

 \bigskip
 
 \bigskip
 
\bigskip

{\bf Acknowledgments.}  
This research was partially supported by the Ministry of Education and Science of the Russian Federation, project
(No 2.3208.2017/PCH), by the Russian Federal Professor program of the Ministry of Education and Science of the Russian Federation,
(project No 1.472.2016/FPM) and by 
 the Academic D.I. Mendeleev Fund Program of the  Tomsk State University (research project NU 8.1.55.2015 L).

\bigskip

\renewcommand{\theequation}{A.\arabic{equation}}
\renewcommand{\thetheorem}{A.\arabic{theorem}}
\renewcommand{\thesubsection}{A.\arabic{subsection}}
\section{Appendix}\label{sec:A}
\setcounter{equation}{0}
\setcounter{theorem}{0}

\subsection{Property of the penalty term}

\begin{lemma}\label{Le.sec:A.1-06-11-01}
For any $n\ge\,1$ and $\lambda \in \Lambda$,
$$ 
P^{0}_n(\lambda) \leq \E_\zs{Q} \Er_n(\lambda)+\frac{\C_\zs{1,Q,n}}{n},
$$
where the coefficient $P^{0}_n(\lambda)$ was defined in \eqref{sec:Pr.2}.
\end{lemma}
\proof
 By the definition of $\Er_n(\lambda)$
  one has
$$ 
\Er_n(\lambda)= \sum_\zs{j=1}^{n} \left((\lambda(j)-1) \theta_\zs{j}+ \frac{\lambda(j)}{n}\xi_\zs{j,n} \right)^2
\,. 
$$                                                                               
In view of Proposition \ref{Pr.sec:OI.2}, this leads to the desired result
$$
\E_\zs{Q}\, \Er_n(\lambda) \ge\, \frac{1}{n}\sum_\zs{j=1}^{n} \lambda^2(j)  \E_\zs{Q}\,\xi^2_\zs{j,n} \ge\, P^{0}_n(\gamma)-\frac{\C_\zs{1,Q,n}}{n}. 
$$
\fdem

\subsection{Properties of the Fourier transform}\label{sec:Four}

\begin{theorem}\label{Th.sec:A.1-00} Cauchy (1825)

Let $U$ be a simply connected open subset of $\bbc,$ let $g : U \to \bbc$ be a holomorphic function, and let $\gamma$ be a rectifiable path in $U$
 whose start point is equal to its end point. Then
$$
\oint_\zs{\gamma}\,g(z)\d z=0\,.
$$
\end{theorem}

\begin{proposition}\label{Pr.sec:A.1-00} 
Let $g : \bbc\to\bbc$ be a holomorphic function in 
$U=\left\{z\in\bbc\,:\,-\beta_\zs{1} < \mbox{Im} z< \beta_\zs{2}\right\}$
for some $\beta_\zs{1}>0$ and $\beta_\zs{2}>0$. Assume that, for any $-\beta_\zs{1}\le t\le 0,$
\begin{equation}\label{sec:A.7-00}
\int_\zs{\bbr}\,\vert g(\theta+i t)\vert\,\d \theta<\infty
\quad\mbox{and}\quad
\lim_\zs{\vert\theta\vert\to\infty}\,g(\theta+i t)\,=0\,.
\end{equation}
Then, for any $x\in\bbr$ and for any $0<\beta<\beta_\zs{1},$
\begin{equation}\label{sec:A.7-0}
\int_\zs{\bbr}\,e^{i\theta x} g(\theta)\,\d \theta
=e^{-\beta x}
\int_\zs{\bbr}\,e^{i\theta x} g(\theta-i\beta )\,\d \theta.
\end{equation}
\end{proposition}
\proof
First note that the conditions of this theorem imply that
$$
\int_\zs{\bbr}\,e^{i\theta x} g(\theta)\,\d \theta
=
\lim_\zs{N\to\infty}\,
\int^{N}_\zs{-N}\,e^{i\theta x} g(\theta)\,\d \theta\,.
$$
We fix now $0<\beta<\beta_\zs{1}$ and we set for any $N\ge 1$ 
\begin{align*}
\gamma&=\{z\in\bbc:-N\le \mbox{Re}z\le N\,,\,\mbox{Im}z=0\}
\cup
\{z\in\bbc:-N\le \mbox{Im}z\le N\,,\,\mbox{Re}z=N\}\\[2mm]
&\cup
\{z\in\bbc:-N\le \mbox{Re}z\le N\,,\,\mbox{Im}z=-\beta\}
\cup
\{z\in\bbc:-\beta\le \mbox{Im}z\le 0\,,\,\mbox{Re}z=-N\}\,.
\end{align*}
Now, in view of the Cauchy theorem, we obtain that
for any $N\ge 1$
\begin{align}\nonumber
\oint_\zs{\gamma}\,&e^{iz x}\,g(z)\d z=
\int^{N}_\zs{-N}\,e^{i\theta x} g(\theta)\,\d \theta
+
\int^{-\beta}_\zs{0}\,e^{i(N+i t)x} g(N+i t)\,\d t
\\[2mm]\label{sec:A.8-00}
&
+
\int^{-N}_\zs{N} e^{i(-i\beta+\theta ) x} g(-i\beta+\theta)\d \theta
+
\int^{0}_\zs{-\beta}e^{i(-N+i t)x} g(-N+i t)\d t
=0\,.
\end{align}
The conditions \eqref{sec:A.7-00} provide that
$$
\lim_\zs{N\to\infty}\,\int^{-\beta}_\zs{0}\,e^{i(N+i t)x} g(N+i t)\,\d t=
\lim_\zs{N\to\infty}\,\int^{0}_\zs{-\beta}\,e^{i(-N+i t)x} g(-N+i t)\,\d t=0\,.
$$
Therefore, letting $N\to\infty$ in \eqref{sec:A.8-00} we obtain \eqref{sec:A.7-0}.
Hence we get the desired result.
\fdem

\bigskip

The following technical lemma is also needed in the present paper.

\begin{lemma}\label{Le.sec:A.2-00-1} 
Let $g : [a,b]\to\bbr$ be a function from $\L_\zs{1}[a,b]$. 
Then, for any fixed $-\infty\le a<b\le +\infty,$
\begin{equation}
\label{sec:A.8-00-1}
\lim_\zs{N\to\infty}\,\int^{b}_\zs{a}\,g(x)\,\sin(Nx)\d x=0
\quad\mbox{and}\quad
\lim_\zs{N\to\infty}\,\int^{b}_\zs{a}\,g(x)\,\cos(Nx)\d x=0
\,.
\end{equation}
\end{lemma}
\proof Let first $-\infty< a<b< +\infty$.
Assume  that $g$ is continuously differentiable, i.e.
 $g\in\C^{1}[a,b]$. Then integrating by parts gives us
 $$
 \int^{b}_\zs{a}\,g(x)\,\sin(Nx)\,\d x
 =\frac{1}{N}\left(
 g(b)\,\sin(Nb)\
 -
 g(a)\,\sin(Na)\
 -
  \int^{b}_\zs{a}\,g^{'}(x)\,\cos(Nx)\,\d x
 \right)\,.
 $$
 So, from this we obtain that
 $$
 \left
 \vert
 \int^{b}_\zs{a}\,g(x)\,\sin(Nx)\,\d x
 \right
 \vert
 \le \frac{\vert g(a)\vert+\vert g(a)\vert+(b-a)\max_\zs{a\le x\le b}\vert g^{'}(x)\vert}{N}
 \,.
 $$
 This implies the first limit in \eqref{sec:A.8-00-1} for this case. The second one is obtained similarly.
Let now $g$ be any absolutely integrated function on $[a,b]$, i.e. $g\in\L_\zs{1}[a,b]$. In this case there exists 
a sequence $g_\zs{n}\in\C^{1}[a,b]$ such that
$$
\lim_\zs{n\to\infty}\int^{b}_\zs{a}\vert g(x)-g_\zs{n}(x)\vert \d x=0\,.
$$
Therefore, taking into account that for any $n\ge 1$
$$
\lim_\zs{N\to\infty}\,\int^{b}_\zs{a}\,g_\zs{n}(x)\,\sin(Nx)\d x=0\,,
$$
we obtain that
$$
\limsup_\zs{n\to\infty}\,\vert \int^{b}_\zs{a}\,g(x)\,\sin(Nx)\d x\vert
\le 
\int^{b}_\zs{a}\vert g(x)-g_\zs{n}(x)\vert \d x\,.
$$
So, letting in this inequality $n\to\infty$ we obtain the first limit in \eqref{sec:A.8-00-1} and, similarly, we obtain the second one.
Let now $b=+\infty$ and $a=-\infty$. In this case we obtain that for any $-\infty<a<b<+\infty$
\begin{align*}
\left\vert 
\int^{+\infty}_\zs{-\infty}\,g(x)\,\sin(Nx)\d x
\right\vert
\le& 
\left\vert 
\int^{+\infty}_\zs{-\infty}\,g(x)\,\sin(Nx)\d x
\right\vert
+
\int^{+\infty}_\zs{b}\,\vert g(x)\,\vert \d x\\
&+
\int^{a}_\zs{-\infty}\,\vert g(x)\,\vert \d x
\,.
\end{align*}
Using here the previous results we obtain that for any $-\infty<a<b<+\infty$
$$
\limsup_\zs{N\to\infty}
\left\vert 
\int^{+\infty}_\zs{-\infty}\,g(x)\,\sin(Nx)\d x
\right\vert
\le 
\int^{+\infty}_\zs{b}\,\vert g(x)\,\vert \d x
+
\int^{a}_\zs{-\infty}\,\vert g(x)\,\vert \d x
\,.
$$
Passing here to limit as $b\to +\infty$ and $a\to-\infty$ we obtain the first limit in
\eqref{sec:A.8-00-1}. Similarly, we can obtain the second one. \fdem \\

Let us now study the inverse Fourier transform. To this end, we need the following local Dini
condition.

$\D$) {\em Assume that, for some fixed $x\in\bbr,$ 
there exist the finite limits
$$
g(x-)=\lim_\zs{z\to x-}g(z)
\quad\mbox{and}\quad
g(x+)=\lim_\zs{z\to x+}g(z)
$$
and 
there exists $\delta=\delta(x)>0$
for which
\begin{equation*}
\int^{\delta}_\zs{0}\,
\frac{
\vert 
g(x+t)+g(x-t)-g(x+)-g(x-)
\vert
}{t}
\d t\,<\,\infty.
\end{equation*}
}

\begin{proposition}\label{Pr.sec:A.2-00} 
Let $g : \bbr\to\bbr$ be a function from $\L_\zs{1}(\bbr)$. If, for some $x\in\bbr,$ this function satisfies
the condition $\D$, then  
\begin{equation}\label{sec:A.8-01}
g(x+)+g(x-)
=
\frac{1}{\pi}
\int_\zs{\bbr}\,e^{-i\theta x} \wh{g}(\theta)\,\d \theta
\,,
\end{equation}
where
$$
\wh{g}(\theta)=\int_\zs{\bbr}\,e^{i\theta t}\,g(t)\,\d t\,.
$$
\end{proposition}

\proof
First, for any fixed $N>0$ we set
$$
J_\zs{N}(x)=\frac{1}{2\pi}\,
\int^{N}_\zs{-N}\,e^{-i\theta x} \wh{g}(\theta)\,\d \theta
=\frac{1}{\pi}\,
\int_\zs{\bbr}\,g(z)\,
\int^{N}_\zs{0}\,\cos(\theta(z- x))\,\d \theta
\d z\,,
$$
i.e.,
$$
J_\zs{N}(x)=
\frac{1}{\pi}\,
\int_\zs{\bbr}\,g(z)\,
\frac{\sin(N(z- x))}{z- x}
\,
\d z
=
\frac{1}{\pi}\,
\int^{\infty}_\zs{0}\,(g(x+t)+g(x-t))\,
\frac{\sin(N t)}{t}
\,
\d t\,.
$$
Taking into account that for any $N>0$
the integral
\begin{equation}\label{sec:A.9-00-1}
\frac{2}{\pi}\,
\int^{\infty}_\zs{0}\,
\frac{\sin(N t)}{t}
\,
\d t
=1
\end{equation}
and denoting 
$B(x)=(g(x+)+g(x-))/2$,
we obtain that
$$
J_\zs{N}(x)-B(x)
=
\frac{1}{\pi}\,
\int^{\infty}_\zs{0}\,
\frac{\omega(x,t)\sin(N t)}{t}
\,
\d t
\quad\mbox{and}\quad
\omega(x,t)=g(x+t)+g(x-t)-2B(x)\,.
$$
Now we represent the last integral as
$$
\int^{\infty}_\zs{0}\,
\frac{\omega(x,t)\sin(N t)}{t}
\,
\d t
=I_\zs{1,N}
+
I_\zs{2,N}
-
2B(x)
I_\zs{3,N}\,,
$$
where
$$
I_\zs{1,N}=\int^{\delta}_\zs{0}\,
\frac{\omega(x,t)}{t}\,\sin(N t)\d t\,,\quad
I_\zs{2,N}=\int^{\infty}_\zs{\delta}\,
G(t)\,\sin(N t)\d t\,,
\quad
I_\zs{3,N}=\int^{\infty}_\zs{\delta}\,
\frac{\sin(N t)}{t}\,\d t
$$
and $G(t)=(g(x+t)+g(x-t))/t$.
Condition $\D$ and Lemma \ref{Le.sec:A.2-00-1} imply directly the convergence $I_\zs{1,N}\to 0$ as $N\to\infty$.  Now note that,
since $g\in\L_\zs{1}(\bbr),$ then the function  $G$ is absolutely integrated. Therefore, in view of  Lemma \ref{Le.sec:A.2-00-1}, 
$I_\zs{2,N}\to 0$ as $N\to\infty$.  As to the last integral we use the property \eqref{sec:A.9-00-1}, i.e.,
the changing  of the variables gives 
$$
I_\zs{3,N}=\int^{\infty}_\zs{\delta N}\,
\frac{\sin t}{t}\,\d t
\to 0
\quad\mbox{as}\quad
N\to\infty\,.
$$
Hence we have the desired result.   \fdem


\medskip

\medskip

\newpage

\end{document}